\documentclass{amsart}
\usepackage{amscd,amssymb}
\usepackage[all]{xy}
\begin{document}

\makeatletter
\def\revddots{\mathinner{\mkern1mu\raise\p@
  \vbox{\kern7\p@\hbox{.}}\mkern2mu
  \raise4\p@\hbox{.}\mkern2mu\raise7\p@\hbox{.}\mkern1mu}}
\makeatother

\newcommand{\Inn}{\operatorname{Inn}\nolimits}
\newcommand{\Der}{\operatorname{Der}\nolimits}
\newcommand{\Hom}{\operatorname{Hom}\nolimits}
\renewcommand{\Im}{\operatorname{Im}\nolimits}
\newcommand{\Ker}{\operatorname{Ker}\nolimits}
\newcommand{\Coker}{\operatorname{Coker}\nolimits}
\newcommand{\rad}{\operatorname{rad}\nolimits}
\newcommand{\rrad}{\mathfrak{r}}
\newcommand{\Ext}{\operatorname{Ext}\nolimits}
\newcommand{\id}{{\operatorname{id}\nolimits}}
\newcommand{\pd}{{\operatorname{pd}\nolimits}}
\newcommand{\End}{\operatorname{End}\nolimits}
\renewcommand{\L}{\Lambda}
\newcommand{\Y}{{\mathcal Y}}
\renewcommand{\P}{{\mathcal P}}
\newcommand{\HH}{\operatorname{HH}\nolimits}
\newcommand{\mo}{\mathfrak{o}}
\renewcommand{\mod}{\operatorname{mod}\nolimits}
\newcommand{\mt}{\mathfrak{t}}
\newcommand{\NonTip}{\operatorname{NonTip}\nolimits}
\newcommand{\typ}{\operatorname{typ}\nolimits}
\newcommand{\kar}{\operatorname{char}\nolimits}

\newtheorem{lem}{Lemma}[section]
\newtheorem{prop}[lem]{Proposition}
\newtheorem{cor}[lem]{Corollary}
\newtheorem{thm}[lem]{Theorem}
\newtheorem{bit}[lem]{}
\theoremstyle{definition}
\newtheorem{defin}[lem]{Definition}
\newtheorem*{remark}{Remark}
\newtheorem{example}[lem]{Example}




\title[Self-injective algebras and Hochschild cohomology]
{Self-injective algebras and the second Hochschild cohomology group}
\author[Al-Kadi]{Deena Al-Kadi}

\address{Deena Al-Kadi\\
Department of Mathematics\\
University of Leicester\\
University Road\\
Leicester, LE1 7RH\\
England}

\begin{abstract}
\sloppy In this paper we study the second Hochschild cohomology group
${\HH}^2(\L)$ of a finite dimensional algebra $\L$. In particular, we
determine ${\HH}^2(\L)$ where $\L$ is a finite dimensional self-injective
algebra of finite representation type over an algebraically closed field $K$
and show that this group is zero for most such $\L$; we give a basis for
${\HH}^2(\L)$ in the few cases where it is not zero.
\end{abstract}

\date{\today}

\maketitle
\section*{Acknowledgements}
This paper is a part of my PhD.\ thesis at the University of
Leicester. I thank Taif University in Saudi Arabia for funding my
PhD.\ research, and Dr.\ Nicole Snashall for her helpful supervision
and valuable suggestions. The author also thanks the referee for
their helpful comments.

\section*{Introduction}
In this paper we study the second Hochschild cohomology group ${\HH}^2(\L)$
of all finite dimensional self-injective algebras $\L$ of finite
representation type over an algebraically closed field $K$.

In general, finite dimensional self-injective algebras of finite
representation type over an algebraically closed field $K$ were shown by
Riedtmann in \cite{R1} to fall into one of the types  $A$, $D$ or $E$,
depending on the tree class of the stable Auslander-Reiten quiver of the
algebra. Riedtmann classified the stable equivalence representatives of
these algebras of type $A$ in \cite{R}; Asashiba then showed that the stable
equivalence classes are exactly the derived equivalence classes for types
$A$, $D$ in \cite[Theorem 2.2]{A2}. In \cite{A},the derived equivalence
class representatives are given explicitly by quiver and relations.

Happel showed in \cite{H} that Hochschild cohomology is invariant under
derived equivalence. So if $A$ and $B$ are derived equivalent then
${\HH}^2(A) \cong {\HH}^2(B)$. Hence to study ${\HH}^2(\L)$ for all finite
dimensional self-injective algebras of finite representation type over an
algebraically closed field $K$, it is enough to study ${\HH}^2(\L)$ for the
representatives of the derived equivalence classes. The algebras of type $A$
fall into two types: Nakayama algebras and M\"obius algebras, and the
Hochschild cohomology of these algebras has already been studied. In
\cite{EH}, Erdmann and Holm give the dimension of the second Hochschild
cohomology group of a Nakayama algebra. In \cite{GS}, Green and Snashall
find the second Hochschild cohomology group for the M\"obius algebras.

The main work of this paper is thus in determining ${\HH}^2(\L)$ for
the finite dimensional self-injective algebras of finite
representation type of types $D$ and $E$. In Section \ref{ch2} we
give a summary of \cite{A} which gives the explicit derived
equivalence representatives we consider. Section \ref{ch3} gives a
short description of the projective resolution of \cite{GS} which we
use to find ${\HH}^2(\L)$. In Section \ref{ch5}, we give a general
theorem, Theorem \ref{motivated thm 1}, which we use to show that
${\HH}^2(\L) = 0$ for most of our algebras. This is motivated by
work in \cite{GS}. The strategy of the theorem is to show that every
element in ${\Hom}(Q^2, \L)$ is a coboundary so that ${\HH}^2(\L) =
0$, where $Q^2$ is the second projective in a minimal projective
resolution of $\L$ as a $\L$,$\L$-bimodule. For all other cases
which are not covered by Theorem \ref{motivated thm 1}, we determine
${\HH}^2(\L)$ by direct calculation, and find a basis for
${\HH}^2(\L)$ in the instances where ${\HH}^2(\L) \neq 0$. The
standard algebras are considered in Sections \ref{sec4} and
\ref{sec5} and the non-standard algebras in Section \ref{sec6}.
Finally Theorem \ref{finalthm} summarises our results and describes
${\HH}^2(\L)$ for all finite dimensional self-injective algebras
$\L$ of finite representation type over an algebraically closed
field. As a consequence, we show that $\dim\HH^2(\L) \neq
\dim\HH^2(\L')$ for a non-standard algebra $\L$ and its standard
form $\L'$, where $\L$ and $\L'$ are of type $(D_{3m},1/3,1)$. This
gives an alternative proof that $\L$ and $\L'$ are not derived
equivalent.

\section{The derived equivalence representatives}\label{ch2}

We give here Asashiba's full classification from \cite{A} and \cite{A2} of
the derived equivalence class representatives of the finite dimensional
self-injective algebras of finite representation type over an algebraically
closed field. These derived equivalence class representatives are listed
according to their type.

From \cite{R1}, the stable Auslander Reiten quiver of a self-injective
algebra $\L$ of finite representation type has the form ${\mathbb Z} \Delta
/ \langle g \rangle,$ where $\Delta$ is a Dynkin graph, $g = \zeta
\tau^{-r}$ such that $r$ is a natural number, $\zeta$ is an automorphism of
the quiver ${\mathbb Z} \Delta$ with a fixed vertex, and $\tau$ is the
Auslander-Reiten translate. Then $\typ(\L) := (\Delta, f, t)$, where $t$
is the order of $\zeta$ and $f := r / m_{\Delta}$ such that $m_{\Delta} = n,
2n - 3, 11, 17$ or 29 as $\Delta = A_n, D_n, E_6, E_7$ or $E_8$,
respectively. We take the following results from \cite{A2}.

\begin{prop} \cite[Theorem 2.2]{A2} \label{propro1}
Given $\L$ a self-injective algebra of finite representation type then the
type $\typ(\L)$ is an element of one of the following sets:

$\{(A_n, s/n, 1)|n, s \in {\mathbb N}\};$

$\{(A_{2p + 1}, s, 2)|p, s \in {\mathbb N}\};$

$\{(D_n, s, 1)|n, s \in {\mathbb N}, n \geq 4\};$

$\{(D_n, s, 2)|n, s \in {\mathbb N}, n \geq 4\};$

$\{(D_4, s, 3)|s \in {\mathbb N}\};$

$\{(D_{3m}, s/3, 1)|m, s \in {\mathbb N}, m \geq 2, 3 \nmid s\};$

$\{(E_n, s, 1)|n = 6, 7, 8, s \in {\mathbb N}\}; \mbox{ and }$

$\{(E_6, s, 2)|s \in {\mathbb N}\}.$
\end{prop}

\begin{thm} \label{thmthm1} \cite[Theorem 2.2]{A2}
Let $\L$ and $\Pi$ be self-injective algebras of finite representation type.

(i) If $\L$ is standard and  $\Pi$ is non-standard then $\L$ and $\Pi$ are
not derived equivalent.

(ii) If $\L$ and $\Pi$ are either both standard or both non-standard then
the following are equivalent:

1) $\L$ and $\Pi$ are derived equivalent;

2) $\L$ and $\Pi$ are stably equivalent;

3) $\typ(\L) = \typ(\Pi)$.
\end{thm}

\vspace{.5cm}

Using these results, \cite{A} gives the derived equivalence
representatives by quiver and relations; these are stated here for
convenience. The derived equivalence representatives of the standard
algebras are given in \ref{bit1}-\ref{bit8}. The non-standard
derived equivalence representatives are given in \ref{bit9}. Recall
from \cite[Theorem 2.2]{A2} that the non-standard derived
equivalence representatives only occur when $\kar K = 2$. Note that
$[j]$ denotes the residue of $j$ modulo $s$ where $s \geq 1$ and we
write paths from left to right (whereas paths are written from right
to left in \cite{A}).

\begin{bit}\label{bit1}
{\it $\L(A_n, s/n, 1)$ with $s, n  \geq 1$.}
\end{bit}

\smallskip

$\L(A_n, s/n, 1)$ with $s, n \geq 1$ is the Nakayama algebra $N_{s, n}$ and
it is given by the quiver $Q(N_{s, n})$:

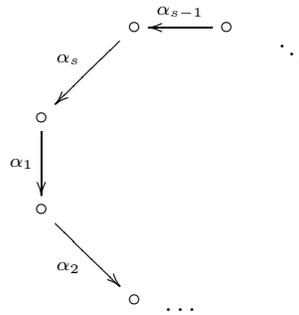
$$\xymatrix{
& \circ\ar[dl]_{\alpha_s} & \circ\ar[l]_{\alpha_{s - 1}}\ar@{}[dr]^\ddots &
\\
\circ\ar[d]_{\alpha_1} & & & \\
\circ\ar[dr]_{\alpha_2} & & &  \\
& \circ\ar@{}[r]_{\displaystyle \ldots} & & \\
}$$

with relations $R(N_{s, n})$:

$\alpha_i \alpha_{i +1} \cdots \alpha_{i + n} = 0$, for all $i \in \{1, 2,
\ldots, s\} = {\mathbb Z} / \langle s \rangle.$

\bigskip

\begin{bit}\label{bit2}
{\it $\L(A_{2p + 1}, s, 2)$ with $s, p  \geq 1$.}
\end{bit}

\smallskip

$\L(A_{2p + 1}, s, 2)$ with $s, p \geq 1$ is the M\"obius algebra $M_{p,s}$
and it is given by the quiver $Q(M_{p, s})$:

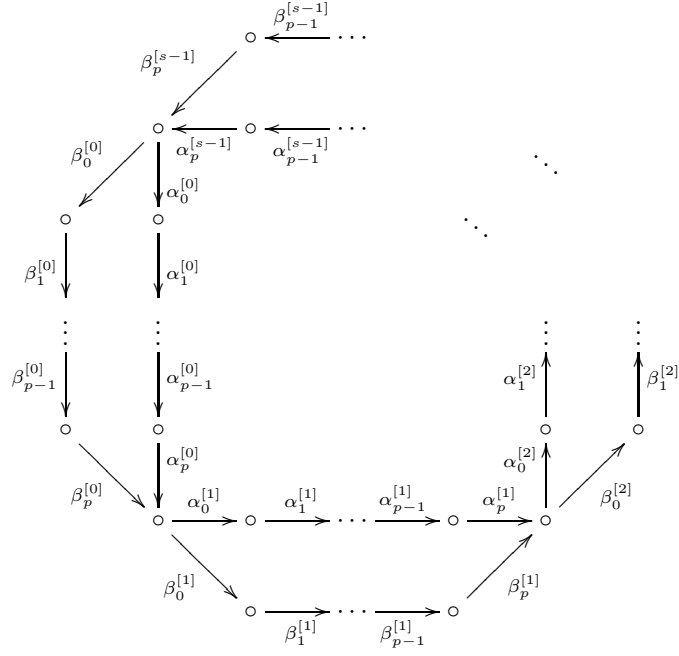
$$\xymatrix{
& & \circ\ar[dl]_{\beta_p^{[s-1]}} & \cdots\ar[l]_{\beta_{p-1}^{[s-1]}} & &
&\\& \circ\ar[dl]_{\beta_0^{[0]}}\ar[d]^>>>{\alpha_0^{[0]}} &
\circ\ar[l]^{\alpha_p^{[s-1]}} & \cdots\ar[l]^{\alpha_{p-1}^{[s-1]}} & &
\ar@{}^\ddots &\\
\circ\ar[d]_{\beta_1^{[0]}} & \circ\ar[d]^{\alpha_1^{[0]}} & & &
\ar@{}^\ddots & &\\
\vdots\ar[d]_{\beta_{p-1}^{[0]}} & \vdots\ar[d]^{\alpha_{p-1}^{[0]}} & & & &
\vdots & \vdots\\
\circ\ar[dr]_{\beta_p^{[0]}} & \circ\ar[d]^<<<{\alpha_p^{[0]}} & & & &
\circ\ar[u]^{\alpha_1^{[2]}} & \circ\ar[u]_{\beta_1^{[2]}}\\
& \circ\ar[r]^{\alpha_0^{[1]}}\ar[dr]_{\beta_0^{[1]}} &
\circ\ar[r]^{\alpha_1^{[1]}} & \cdots\ar[r]^{\alpha_{p - 1}^{[1]}} &
\circ\ar[r]^{\alpha_p^{[1]}} &
\circ\ar[u]^>>>{\alpha_0^{[2]}}\ar[ur]_{\beta_0^{[2]}} & &\\
& & \circ\ar[r]_{\beta_1^{[1]}} & \cdots\ar[r]_{\beta_{p-1}^{[1]}} &
\circ\ar[ur]_{\beta_p^{[1]}} & &\\
}$$

with relations $R(M_{p,s})$:

(i) $\alpha_0^{[i]} \cdots \alpha_p^{[i]} = \beta_0^{[i]} \cdots
\beta_p^{[i]}$, for all $i \in \{0, \ldots, s - 1\}$,

(ii)  for all $i \in \{0, \ldots, s - 2\}$,
$$\alpha_p^{[i]} \beta_0^{[i + 1]} = 0,  \hspace{1cm} \beta_p^{[i]}
\alpha_0^{[i + 1]} = 0,$$
$$\alpha_p^{[s - 1]} \alpha_0^{[0]} = 0, \hspace{1cm} \beta_p^{[s - 1]}
\beta_0^{[0]} = 0,$$

(iii) paths of length $p + 2$ are equal to 0.

\bigskip

\begin{bit}\label{bit3}
{\it $\L(D_n, s, 1)$ with $n \geq 4, s \geq 1$.}
\end{bit}

\smallskip

$\L(D_n, s, 1)$ with $n \geq 4, s \geq 1$ is given by the quiver $Q(D_n,
s)$:

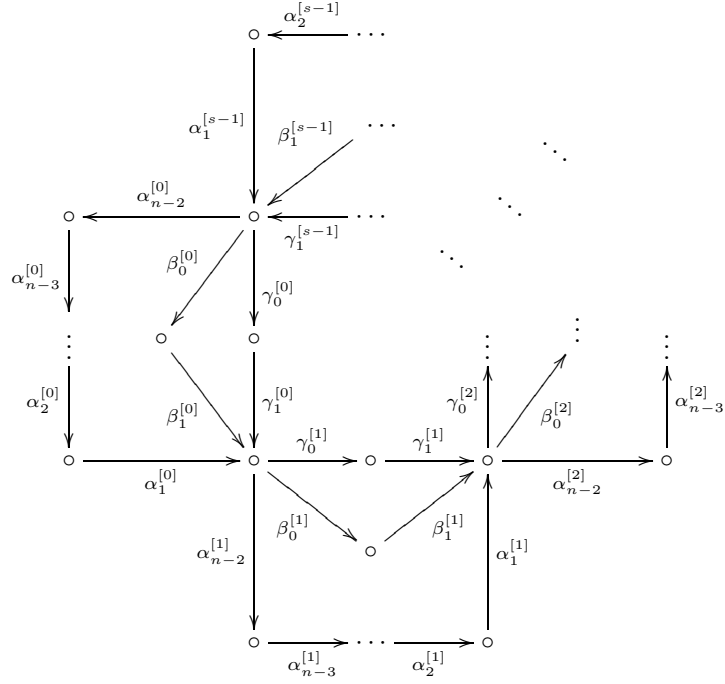
$$\xymatrix{
& & \circ\ar[dd]_{\alpha_1^{[s-1]}} & \cdots\ar[l]_{\alpha_2^{[s-1]}} & &
&\\
& & & \mbox{$\ \ \ \cdots$}\ar[dl]_<<<<{\beta_1^{[s-1]}} &
\ar@{}[dr]^\ddots  & &\\
\circ\ar[d]_{\alpha_{n-3}^{[0]}} & &
\circ\ar[ll]_{\alpha_{n-2}^{[0]}}\ar[dl]_{\beta_0^{[0]}}\ar[d]^>>>>>{\gamma_0^{[0]}}
& \cdots\ar[l]^{\gamma_1^{[s-1]}} \ar@{}[dr]^\ddots & & \ar@{}[ul]^\ddots  &
\\
\vdots\ar[d]_{\alpha_2^{[0]}} & \circ\ar[dr]_{\beta_1^{[0]}} &
\circ\ar[d]^{\gamma_1^{[0]}} & & \vdots & \raisebox{3ex}{\vdots} &  \vdots\\
\circ\ar[rr]_{\alpha_1^{[0]}} & &
\circ\ar[r]^{\gamma_0^{[1]}}\ar[dr]_{\beta_0^{[1]}}\ar[dd]_{\alpha_{n-2}^{[1]}}
& \circ\ar[r]^{\gamma_1^{[1]}} &
\circ\ar[u]^{\gamma_0^{[2]}}\ar[ur]_{\beta_0^{[2]}}\ar[rr]_{\alpha_{n-2}^{[2]}}
& & \circ \ar[u]_{\alpha_{n-3}^{[2]}}\\
& & & \circ\ar[ur]_{\beta_1^{[1]}} & & &\\
& & \circ\ar[r]_{\alpha_{n-3}^{[1]}} & \cdots\ar[r]_{\alpha_2^{[1]}} &
\circ\ar[uu]_{\alpha_1^{[1]}} & &\\
}$$

with relations $R(D_n, s, 1)$:

(i) $\alpha_{n - 2}^{[i]}\alpha_{n - 3}^{[i]} \cdots
\alpha_{2}^{[i]}\alpha_{1}^{[i]} = \beta_{0}^{[i]} \beta_{1}^{[i]} =
\gamma_{0}^{[i]} \gamma_{1}^{[i]}$, for all $i \in \{0, \ldots, s-1\} =
{\mathbb Z} / \langle s \rangle$,

(ii) for all $i \in \{0, \ldots, s-1\} =  {\mathbb Z} / \langle s \rangle$,
$$\alpha_1^{[i]} \beta_0^{[i + 1]}= 0,  \hspace{1cm}
\alpha_1^{[i]}\gamma_0^{[i + 1]}  = 0,$$
$$\beta_1^{[i]}\alpha_{n - 2}^{[i + 1]}  = 0,  \hspace{1cm}
\gamma_1^{[i]}\alpha_{n - 2}^{[i + 1]} = 0,$$
$$\beta_1^{[i]}\gamma_0^{[i + 1]} = 0,  \hspace{1cm}
\gamma_1^{[i]}\beta_0^{[i + 1]} = 0,$$

(iii) for all $i \in \{0, \ldots, s-1\} =  {\mathbb Z} / \langle s \rangle$
and for all $j \in \{1, \ldots, n - 2\} = {\mathbb Z} / \langle n - 2
\rangle$,
$$\alpha_j^{[i]} \ldots \alpha_{j - n + 2}^{[i + 1]} = 0,$$
$$\beta_0^{[i]} \beta_1^{[i]} \beta_0^{[i +1]} = 0, \hspace{1cm}
\gamma_0^{[i]} \gamma_1^{[i]} \gamma_0^{[i +1]} = 0,$$
$$\beta_1^{[i]} \beta_0^{[i + 1]} \beta_1^{[i +1]} = 0, \hspace{1cm}
\gamma_1^{[i]} \gamma_0^{[i + 1]} \gamma_1^{[i +1]} = 0.$$

The set of relations (iii) means that ``$\alpha$-paths'' of length $n - 1$
are equal to 0, ``$\beta$-paths'' of length 3 are equal to 0 and
``$\gamma$-paths'' of length 3 are equal to 0.


\begin{bit}\label{bit4}
{\it $\L(D_n, s, 2)$ with $n \geq 4, s \geq 1$.}
\end{bit}

\smallskip

$\L(D_n, s, 2)$ with $n \geq 4, s \geq 1$ is given by the quiver $Q(D_n, s)$
above with relations $R(D_n, s, 2)$:

(i) $\alpha_{n - 2}^{[i]}\alpha_{n - 3}^{[i]} \cdots
\alpha_{2}^{[i]}\alpha_{1}^{[i]} = \beta_{0}^{[i]} \beta_{1}^{[i]} =
\gamma_{0}^{[i]} \gamma_{1}^{[i]}$, for all $i \in \{0, \ldots, s-1\} =
{\mathbb Z} / \langle s \rangle$,

(ii) for all $i \in \{0, \ldots, s-1\} =  {\mathbb Z} / \langle s \rangle$,
$$ \alpha_{1}^{[i]} \beta_{0}^{[i + 1]}= 0,  \hspace{1cm}
\alpha_{1}^{[i]}\gamma_{0}^{[i + 1]}  = 0,$$
$$\beta_{1}^{[i]}\alpha_{n - 2}^{[i + 1]}  = 0,  \hspace{1cm}
\gamma_{1}^{[i]}\alpha_{n - 2}^{[i + 1]} = 0,$$

and for all $i \in \{0, \ldots, s - 2\}$,
$$ \beta_{1}^{[i]}\gamma_{0}^{[i + 1]} = 0,  \hspace{1cm}
\gamma_{1}^{[i]}\beta_{0}^{[i + 1]} = 0,$$
$$ \beta_{1}^{[s - 1]} \beta_{0}^{[0]} = 0,  \hspace{1cm} \gamma_{1}^{[s -
1]}\gamma_{0}^{[0]} = 0,$$

(iii) ``$\alpha$-paths'' of length $n - 1$ are equal to $0$, and for all $i
\in \{0, \ldots, s - 2\},$
$$\beta_{0}^{[i]} \beta_{1}^{[i]} \beta_{0}^{[i +1]} = 0, \hspace{1cm}
\gamma_{0}^{[i]} \gamma_{1}^{[i]} \gamma_{0}^{[i +1]} = 0,$$
$$\beta_{1}^{[i]} \beta_{0}^{[i + 1]} \beta_{1}^{[i +1]} = 0, \hspace{1cm}
\gamma_{1}^{[i]} \gamma_{0}^{[i + 1]} \gamma_{1}^{[i +1]} = 0 \mbox{ and }$$
$$\beta_{0}^{[s - 1]} \beta_{1}^{[s - 1]} \gamma_{0}^{[0]} = 0, \hspace{1cm}
\gamma_{0}^{[s - 1]} \gamma_{1}^{[s - 1]} \beta_{0}^{[0]} = 0,$$
$$\beta_{1}^{[s - 1]} \gamma_{0}^{[0]} \gamma_{1}^{[0]} = 0, \hspace{1cm}
\gamma_{1}^{[s - 1]} \beta_{0}^{[0]} \beta_{1}^{[0]} = 0.$$

\begin{bit}\label{bit5}
{\it $\L(D_4, s, 3)$ with $s \geq 1$.}
\end{bit}

$\L(D_4, s, 3)$ with $s \geq 1$ is given by the quiver $Q(D_4, s)$ above
with relations $R(D_4, s, 3)$:

(i) $\alpha_{0}^{[i]}\alpha_{1}^{[i]} = \beta_{0}^{[i]} \beta_{1}^{[i]} =
\gamma_{0}^{[i]} \gamma_{1}^{[i]}$, for all $i \in \{0, \ldots, s-1\} =
{\mathbb Z} / \langle s \rangle$,

(ii) for all $i \in \{0, \ldots, s-2\}$,
$$ \alpha_{1}^{[i]} \beta_{0}^{[i + 1]}= 0,  \hspace{1cm}
\alpha_{1}^{[i]}\gamma_{0}^{[i + 1]} = 0,$$
$$\beta_{1}^{[i]}\alpha_0^{[i + 1]}  = 0,  \hspace{1cm}
\gamma_{1}^{[i]}\alpha_0^{[i + 1]} = 0,$$
$$ \beta_{1}^{[i]}\gamma_{0}^{[i + 1]} = 0,  \hspace{1cm}
\gamma_{1}^{[i]}\beta_{0}^{[i + 1]} = 0,$$ and

$$\alpha_1^{[s-1]} \alpha_0^{[0]} = 0, \hspace{1cm} \alpha_1^{[s-1]}
\gamma_0^{[0]} = 0,$$
$$\beta_{1}^{[s-1]} \alpha_0^{[0]} = 0, \hspace{1cm}
\beta_1^{[s-1]} \beta_0^{[0]} = 0,$$
$$\gamma_1^{[s-1]} \beta_{0}^{[0]} = 0, \hspace{1cm}
\gamma_{1}^{[s-1]} \gamma_{0}^{[0]} = 0,$$

(iii) paths of length 3 are equal to 0.


\begin{bit}\label{bit6}
{\it $\L(D_{3m}, s/3, 1)$ with $m \geq 2 \mbox{ and } 3 \nmid s \geq 1$.}
\end{bit}

$\L(D_{3m}, s/3, 1)$ with $m \geq 2 \mbox{ and } 3 \nmid s \geq 1$ is given
by the quiver $Q(D_{3m}, s/3)$:

$$
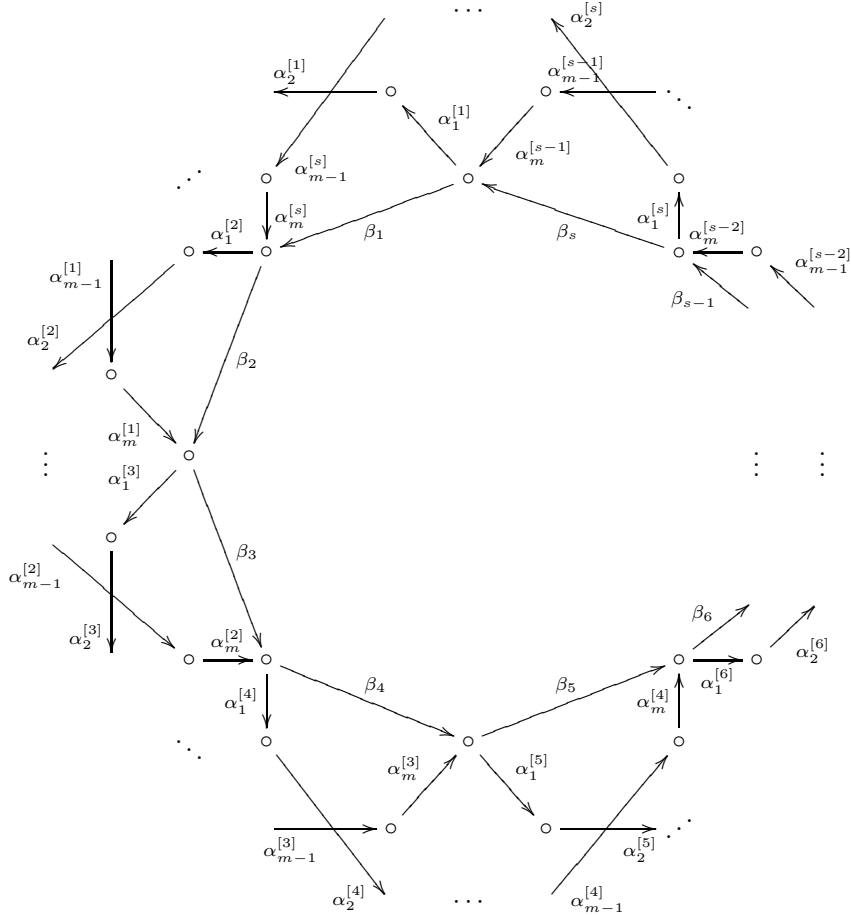
\xymatrix@R=15pt@C=15pt{
& & & & &  \ar[ddll]^>>>{\alpha_{m-1}^{[s]}} & \cdots & & & & & & \\
& & & & &  \circ\ar[ll]_>>>{\alpha_2^{[1]}} & &
\circ\ar[dl]^{\alpha_m^{[s-1]}} & &  \ar[ll]_>>>{\alpha_{m-1}^{[s-1]}}
\ddots & & & \\
& & \revddots & \circ\ar[d]^{\alpha_m^{[s]}} &  & & \circ\ar[dlll]^{\beta_1}
\ar[ul]_{\alpha_1^{[1]}} & & & \circ\ar[uull]_>>>{\alpha_2^{[s]}} & & & \\
& \ar[dd]_<<<{\alpha_{m-1}^{[1]}} & \circ\ar[ddll]_>>>{\alpha_2^{[2]}} &
\circ\ar[l]_{\alpha_1^{[2]}} \ar[dddl]^{\beta_2} & & & & & &
\circ\ar[ulll]^{\beta_s} \ar[u]^{\alpha_1^{[s]}} &
\circ\ar[l]_{\alpha_m^{[s-2]}} & & \\
& & & & & & & & & &  \ar[ul]^{\beta_{s-1}} &  \ar[ul]_{\alpha_{m-1}^{[s-2]}}
& \\
& \circ\ar[dr]_{\alpha_m^{[1]}} & & & & & & & & & & & \\
\vdots & & \circ\ar[dl]_{\alpha_1^{[3]}} \ar[dddr]^{\beta_3} & & & & & & & &
\vdots & \vdots & \\
\ar[ddrr]_<<<{\alpha_{m-1}^{[2]}} & \circ\ar[dd]_>>>{\alpha_2^{[3]}} & & & &
& & & & & & & \\
& & & & & & & & & & &  &  \\
& & \circ\ar[r]^{\alpha_m^{[2]}} & \circ\ar[drrr]^{\beta_4}
\ar[d]_{\alpha_1^{[4]}} & & & & & & \circ\ar[ur]^{\beta_6}
\ar[r]_{\alpha_1^{[6]}} & \circ\ar[ur]_{\alpha_2^{[6]}} & & \\
& & \ddots  & \circ\ar[ddrr]_>>>{\alpha_2^{[4]}} & & &
\circ\ar[urrr]^{\beta_5} \ar[dr]^{\alpha_1^{[5]}} & & &
\circ\ar[u]^{\alpha_m^{[4]}} & & & \\
& & & \ar[rr]_<<<{\alpha_{m-1}^{[3]}} & & \circ\ar[ur]^{\alpha_m^{[3]}} & &
\circ\ar[rr]_>>>{\alpha_2^{[5]}}  & & \revddots  & & & \\
& & & & & & \cdots &  \ar[uurr]_<<<{\alpha_{m-1}^{[4]}} & & & & & \\
}$$

and for $s = 1$, $Q(D_{3m}, 1/3)$:

$$\xymatrix{
& m \ar[dl]_{\alpha_m} & m-1 \ar[l]_{\alpha_{m - 1}} \ar@{}[dr]^\ddots & \\
1 \ar@(u,l)[]_{\beta}\ar[d]_{\alpha_1} & & & \\
2 \ar[dr]_{\alpha_2} & & & \\
& 3 \ar@{}[r]_{\displaystyle \ldots} & & \\
}$$

with relations $R(D_{3m}, s/3, 1)$:

(i) $\alpha_1^{[i]} \alpha_2^{[i]} \cdots \alpha_m^{[i]} = \beta_i
\beta_{i+1}$, for all $i \in \{1, \ldots, s\} = {\mathbb Z} / \langle s
\rangle$,

(ii) $\alpha_m^{[i]} \alpha_1^{[i + 2]} = 0$, for all $i \in \{1, \ldots,
s\} =  {\mathbb Z} / \langle s \rangle$,

(iii) $\alpha_j^{[i]} \cdots \alpha_m^{[i]} \beta_{i+2} \alpha_1^{[i+3]}
\cdots \alpha_j^{[i+3]} = 0$ , for all $i \in \{1, \ldots, s\} =  {\mathbb
Z} / \langle s \rangle$ and for all $j \in \{1, \ldots, m\}$ (i.e. paths of
length $m + 2$ are equal to 0).

\smallskip

In the case $s=1$, the relations $R(D_{3m}, 1/3, 1)$ are:

(i) $\alpha_1 \alpha_2 \cdots \alpha_m = \beta^2$,

(ii) $\alpha_m \alpha_1 = 0$,

(iii) $\alpha_j \cdots \alpha_m \beta \alpha_1 \cdots \alpha_j = 0$ for $j =
2, \ldots, m-1$.

\bigskip

\begin{bit}\label{bit7}
{\it $\L(E_n, s, 1)$ with $n  \in \{6, 7, 8\} \mbox{ and } s \geq 1$.}
\end{bit}

\smallskip

$\L(E_n, s, 1)$ is given by the quiver $Q(E_n, s)$:

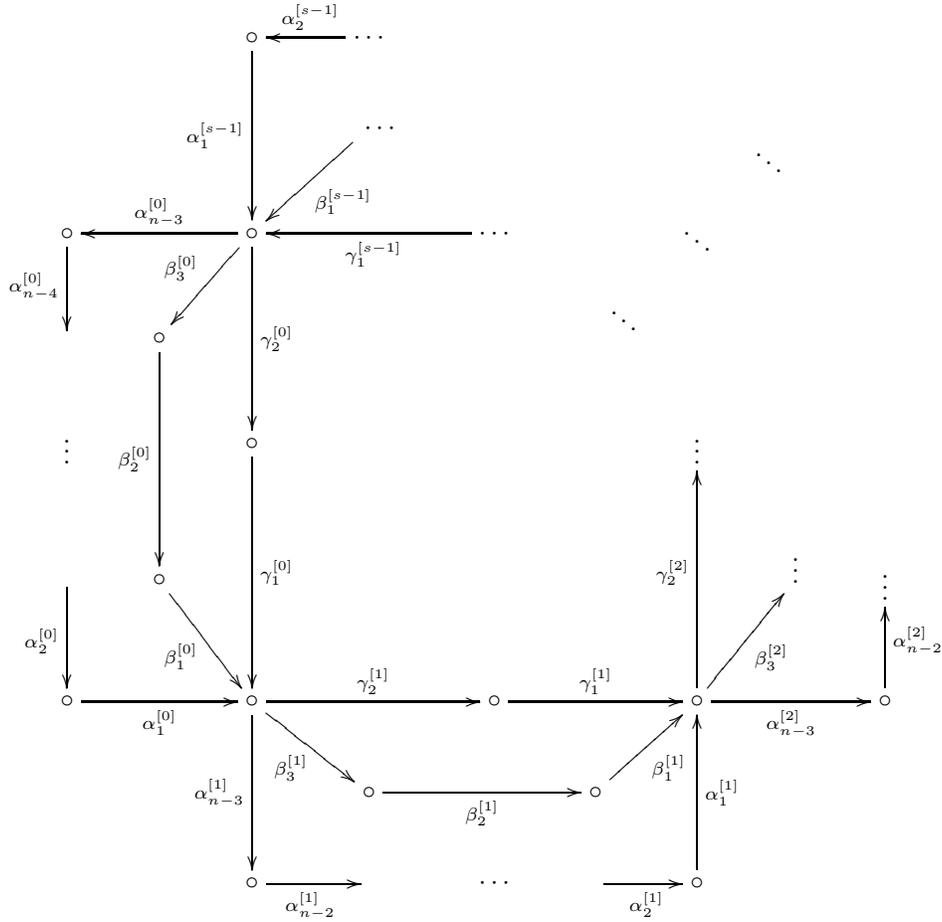
$$\xymatrix{
& & \circ\ar[dd]_{\alpha_1^{[s-1]}} & \ar[l]_{\alpha_2^{[s-1]}} \cdots & & &
& & \\
& & & \mbox{$\ \ \ \cdots$}\ar[dl]^{\beta_1^{[s-1]}} & & & \ar@{}[dr]^\ddots
& & \\
\circ\ar[d]_{\alpha_{n-4}^{[0]}} & &
\circ\ar[ll]_{\alpha_{n-3}^{[0]}}\ar[dl]_{\beta_3^{[0]}}\ar[dd]^{\gamma_2^{[0]}}
& & \cdots\ar[ll]^{\gamma_1^{[s-1]}} & & \ddots & & \\
& \circ\ar[dd]_{\beta_2^{[0]}} & & & & & \ar@{}[ul]^\ddots  & & \\
\vdots & & \circ\ar[dd]^{\gamma_1^{[0]}} &  &  & & \vdots & & \\
\ar[d]_{\alpha_2^{[0]}} & \circ\ar[dr]_{\beta_1^{[0]}} & & & & & &
\raisebox{3ex}{\vdots} & \vdots \\
\circ\ar[rr]_{\alpha_1^{[0]}} & &
\circ\ar[rr]^{\gamma_2^{[1]}}\ar[dr]_{\beta_3^{[1]}}\ar[dd]_{\alpha_{n-3}^{[1]}}
& & \circ\ar[rr]^{\gamma_1^{[1]}} & &
\circ\ar[uu]^{\gamma_2^{[2]}}\ar[ur]_{\beta_3^{[2]}}\ar[rr]_{\alpha_{n-3}^{[2]}}
& & \circ \ar[u]_{\alpha_{n-2}^{[2]}} \\
& & & \circ\ar[rr]_{\beta_2^{[1]}} & & \circ\ar[ur]_{\beta_1^{[1]}} & & & \\
& & \circ\ar[r]_{\alpha_{n-2}^{[1]}} & & \cdots & \ar[r]_{\alpha_2^{[1]}} &
\circ\ar[uu]_{\alpha_1^{[1]}} & & \\
}$$

with relations $R(E_n, s, 1):$

(i) $\alpha_{n-3}^{[i]} \cdots \alpha_2^{[i]} \alpha_1^{[i]} = \beta_3^{[i]}
\beta_2^{[i]} \beta_1^{[i]} = \gamma_2^{[i]} \gamma_1^{[i]}$, for all $i \in
\{0, \ldots, s - 1\},$

(ii) for all $i \in \{0, \ldots, s - 1\} = {\mathbb Z} / \langle s \rangle$,
$$\alpha_{1}^{[i]} \beta_3^{[i + 1]}= 0,  \hspace{1cm}
\alpha_{1}^{[i]}\gamma_2^{[i + 1]}  = 0,$$
$$\beta_{1}^{[i]}\alpha_{n - 3}^{[i + 1]}  = 0,  \hspace{1cm}
\gamma_{1}^{[i]}\alpha_{n - 3}^{[i + 1]} = 0,$$
$$\beta_{1}^{[i]}\gamma_2^{[i + 1]} = 0,  \hspace{1cm}
\gamma_{1}^{[i]}\beta_3^{[i + 1]} = 0,$$

(iii) ``$\alpha$-paths'' of length $n - 2$ are equal to 0, ``$\beta$-paths''
of length 4 are equal to 0 and ``$\gamma$-paths'' of length 3 are equal to
0.


\begin{bit}\label{bit8}
{\it $\L(E_6, s, 2)$ with $s \geq 1$.}
\end{bit}

\smallskip

$\L(E_6, s, 2)$ is given by the quiver $Q(E_6, s)$ above with relations
$R(E_6, s, 2):$

(i) $\alpha_3^{[i]} \alpha_2^{[i]} \alpha_1^{[i]} = \beta_3^{[i]}
\beta_2^{[i]} \beta_1^{[i]} = \gamma_2^{[i]} \gamma_1^{[i]}$, for all $i \in
\{0, \ldots, s - 1\},$

(ii) for all $i \in \{0, \ldots, s - 1\} = {\mathbb Z} / \langle s \rangle$,
$$\gamma_{1}^{[i]}\alpha_3^{[i + 1]} = 0, \hspace{1cm} \gamma_1^{[i]}
\beta_3^{[i + 1]} = 0,$$
$$\alpha_1^{[i]} \gamma_2^{[i+1]} = 0, \hspace{1cm}
\beta_{1}^{[i]}\gamma_2^{[i + 1]} = 0,$$

and for all $i \in \{0, \ldots, s - 2\},$
$$\alpha_1^{[i]} \beta_3^{[i+1]} = 0, \hspace{1cm} \beta_1^{[i]}
\alpha_3^{[i+1]} = 0,$$
$$\alpha_1^{[s-1]} \alpha_3^{[0]} = 0, \hspace{1cm} \beta_1^{[s-1]}
\beta_3^{[0]} = 0,$$

(iii) ``$\gamma$-paths'' of length 3 are equal to 0 and for all $i \in \{0,
\ldots, s - 2\}$ and for all $j \in \{1, 2, 3\} = {\mathbb Z} / \langle 3
\rangle$,
$$\alpha_j^{[i]} \cdots \alpha_{j-3}^{[i+1]} = 0, \hspace{1cm}
\beta_j^{[i]} \cdots \beta_{j -3}^{[i+1]} = 0,$$
$$\alpha_j^{[s-1]} \cdots \alpha_1^{[s-1]} \beta_3^{[0]} \cdots
\beta_{j-3}^{[0]} = 0, \hspace{1cm} \beta_j^{[s-1]} \cdots \beta_1^{[s-1]}
\alpha_3^{[0]} \cdots \alpha_{j-3}^{[0]} = 0.$$

\vspace{.5cm}

Thus we have listed all the derived equivalence representatives of the
standard algebras. The derived equivalence representatives of the
non-standard algebras are given next.

\bigskip

\begin{bit}\label{bit9}
{\it $\L(m)$ with $m \geq 2$.}
\end{bit}
In this case $\kar K = 2$ by \cite[Theorem 2.2]{A2}. The
non-standard algebra $\L(m)$ for each $m \geq 2$ is given by the
quiver ${\mathcal Q}(D_{3m}, 1/3)$:

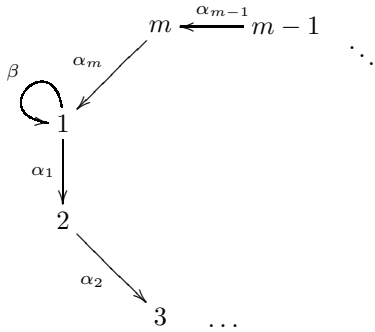
\vspace{1cm}
$$\xymatrix{
& m \ar[dl]_{\alpha_m} & m-1 \ar[l]_{\alpha_{m - 1}} \ar@{}[dr]^\ddots & \\
1 \ar@(u,l)[]_{\beta}\ar[d]_{\alpha_1} & & & \\
2 \ar[dr]_{\alpha_2} & & & \\
& 3 \ar@{}[r]_{\displaystyle \ldots} & & \\
}$$
\vspace{1.5cm}

with relations $R(m)$:

(i) $\alpha_1 \alpha_2 \cdots \alpha_m = \beta^2,$

(ii) $\alpha_m \alpha_1 =  \alpha_m \beta \alpha_1$,

(iii) $\alpha_i \alpha_{i+1} \cdots \alpha_i  = 0$ , for all $i \in \{1,
\ldots, m\} =  {\mathbb Z} / \langle m \rangle$ (i.e. ``$\alpha$''-paths of
length $m + 1$ are equal to 0).

\section{Projective resolutions}\label{ch3}

To find the Hochschild cohomology groups for any finite dimensional algebra
$\L$, a projective resolution of $\L$ as a $\L, \L$-bimodule is needed. In
this section we look at the projective resolutions of \cite{GS} and \cite
{GSZ} in order to describe the second Hochschild cohomology group. Let $\L =
K {\mathcal Q}$/$I$ where ${\mathcal Q}$ is a quiver, and $I$ is an
admissible ideal of $K {\mathcal Q}$. Fix a minimal set $f^2$ of generators
for the ideal $I$. Let $x$ be one of the minimal relations. Then $x =
\sum_{j=1}^{r} c_ja_{1j} \cdots a_{kj} \cdots a_{s_jj}$, that is, $x$ is a
linear combination of paths $a_{1j} \cdots a_{kj} \cdots a_{s_jj}$ for $j =
1, \ldots, r$ and $c_j \in K$ and there are unique vertices $v$ and $w$ such
that each path $a_{1j} \cdots a_{kj} \cdots a_{s_j j}$ starts at $v$ and
ends at $w$ for all $j$. We write $\mo(x) = v$ and $\mt(x) = w.$ Similarly
$\mo(a)$ is the origin of the arrow $a$ and $\mt(a)$ is the end of $a$.

In \cite[Theorem 2.9]{GS}, a minimal projective resolution of $\L$ as a $\L,
\L$-bimodule is given which begins:
$$\cdots \rightarrow Q^3 \stackrel{A_3}{\rightarrow} Q^2
\stackrel{A_2}{\rightarrow} Q^1 \stackrel{A_1}{\rightarrow} Q^0
\stackrel{g}{\rightarrow} \L \rightarrow 0,$$
where the projective $\L, \L$-bimodules $Q^0, Q^1, Q^2$ are given by
$$Q^0 = \bigoplus_{v, vertex} \L v \otimes v\L,$$
$$Q^1 = \bigoplus_{a, arrow} \L \mo(a) \otimes \mt(a)\L, \mbox { and }$$

$$Q^2 = \bigoplus_{x \in f^2} \L \mo(x) \otimes \mt(x) \L.$$

The maps $g, A_1$, $A_2$ and $A_3$ are all $\L, \L$-bimodule homomorphisms.
The map $g: Q^0 \rightarrow \L$ is the multiplication map so is given by  $v
\otimes v \mapsto v$. The map $A_1: Q^1\rightarrow Q^0$ is given by $\mo(a)
\otimes \mt(a) \mapsto \mo(a) \otimes \mo(a) a - a \mt(a) \otimes \mt(a)$
for each arrow $a$.

With the notation for $x \in f^2$ given above, the map $A_2: Q^2 \rightarrow
Q^1$ is given by $\mo(x) \otimes \mt(x) \mapsto
\sum_{j=1}^{r}c_j(\sum_{k=1}^{s_j} a_{1j} \cdots a_{(k-1)j} \otimes
a_{(k+1)j} \cdots a_{s_j j})$, where $a_{1j} \cdots a_{(k-1)j} \otimes
a_{(k+1)j} \cdots a_{s_j j} \in \L \mo(a_{kj}) \otimes \mt(a_{kj})\L$.

In order to find the projective $\L$,$\L$-bimodule $Q^3$ and the map $A_3$
in the $\L, \L$-bimodule resolution of $\L$ in \cite{GS}, Green and Snashall
start by finding a projective resolution of $\L/\rrad$ as a right
$\L$-module, where $\rrad = J(\L)$ is the Jacobson radical of $\L$, using
the notation and procedure of the paper \cite{GSZ}. In \cite{GSZ}, Green,
Solberg and Zacharia show that there are sets $f^n$, $n \geq 3$, and uniform
elements $y \in f^n$ such that $y = \sum_{x \in f^{n-1}} x r_x = \sum_{z \in
f^{n-2}} z s_z$ for unique elements $r_x, s_z \in K{\mathcal Q}$ with
special properties related to a minimal projective $\L$-resolution of
$\L/\rrad$ considered as a right $\L$-module. In particular, for $y \in f^3$
we have $y \in \coprod f^2 K{\mathcal Q} \cap \coprod f^1 I$ and $y$ may be
written $y = \sum f^2_i p_i = \sum q_i f^2_i r_i$ with $p_i, q_i, r_i \in
K{\mathcal Q}$ and $p_i, q_i$ in the ideal generated by the arrows of
$K{\mathcal Q}$ such that the elements $p_i$ are unique. Recall that an
element $y \in K{\mathcal Q}$ is uniform if there are vertices $v, w$ such
that $y = v y = y w.$ We write $\mo(y) = v$ and $\mt(y) = w$.

Then \cite{GS} gives that $Q^3 = \coprod_{y \in f^3} \L \mo(y) \otimes
\mt(y) \L$ and describes the map $A_3$. For $y \in f^3$ in the notation
above, the component of $A_3 (\mo(y) \otimes \mt(y))$ in the summand $\L
\mo(f_i^2) \otimes \mt(f_i^2) \L$ of $Q^2$ is $\Sigma (\mo(y) \otimes p_i -
q_i \otimes r_i).$

Thus we can describe the part of the minimal projective $\L, \L$-bimodule
resolution of $\L$:
$$Q^3 \stackrel{A_3}{\rightarrow} Q^2 \stackrel{A_2}{\rightarrow} Q^1
\stackrel{A_1}{\rightarrow} Q^0 \stackrel{g}{\rightarrow} \L \rightarrow
0.$$
Applying ${\Hom}(-, \L)$ to this resolution gives us the complex
$$0 \rightarrow {\Hom}(Q^0, \L) \stackrel{d_1}{\rightarrow} {\Hom}(Q^1, \L)
\stackrel{d_2}{\rightarrow} {\Hom}(Q^2, \L) \stackrel{d_3}{\rightarrow}
{\Hom}(Q^3, \L)$$
where $d_i$ is the map induced from $A_i$ for $i = 1, 2, 3$. Then
${\HH}^2(\L) = {\Ker}\,d_3/{\Im}\,d_2.$

\bigskip

Throughout, all tensor products are tensor products over $K$, and we write
$\otimes$ for $\otimes_K$. When considering an element of the projective
$\L, \L$-bimodule $Q^1 = \bigoplus_{a, arrow} \L \mo(a) \otimes \mt(a) \L$
it is important to keep track of the individual summands of $Q^1$. So to
avoid confusion we usually denote an element in the summand $\L \mo(a)
\otimes \mt(a) \L$ by $\lambda \otimes_a \lambda'$ using the subscript `$a$'
to remind us in which summand this element lies. Similarly, an element
$\lambda \otimes_{f^2_i} \lambda'$ lies in the summand $\L \mo(f^2_i)
\otimes \mt(f^2_i) \L$ of $Q^2$ and an element $\lambda \otimes_{f^3_i}
\lambda'$ lies in the summand $\L \mo(f^3_i) \otimes \mt(f^3_i) \L$ of
$Q^3$. We keep this notation for the rest of the paper.

\bigskip

Now we are ready to compute ${\HH}^2(\L)$ for the derived equivalence
representatives of the finite dimensional self-injective algebras of finite
representation type over an algebraically closed field.

First we recall that the algebras of type $(A_n, s/n, 1)$ and $(A_{2p+1}, s,
2)$ have been considered in \cite{EH} and \cite{GS} respectively.

\begin{thm} \cite[Theorem 4.2]{GS}
For the M\"obius algebra $M_{p, s}$ we have ${\HH}^2(M_{p, s}) = 0$ except
when $p = 1$ and $s = 1.$
\end{thm}

It is well-known that if $p = 1$ and $s = 1$ then $M_{p, s}$ is the
preprojective algebra of type $A_3$. In \cite{ES}, a basis for the
Hochschild cohomology groups of the preprojective algebras of type $A_n$ is
given.

\begin{prop}\cite[7.2.1]{ES}
For the M\"obius algebra $M_{p, s}$ with $p = 1$ and $s = 1$ we have
$\dim\,{\HH}^2(M_{p, s}) = 1.$
\end{prop}

In \cite{EH}, the dimension of ${\HH}^{2j}(\L)$ is given for a
self-injective Nakayama algebra for all $j \geq 1.$ In particular this gives
us ${\HH}^2(\L)$ when $j = 1.$ The self-injective Nakayama algebra
$\L(A_n,s/n, 1)$ of \cite{A} is the algebra $B_s^{n+1}$ of \cite{EH}. Write
$n+1 = ms + r$ where $0 \leq r < s$. From \cite{EH}, with $j = 1$, we have
the following result.

\begin{prop} \cite[Proposition 4.4]{EH}
For $\L = \L(A_n, s/n, 1)$, and with the above notation we have
$\dim\,{\HH}^2(\L) = m.$
\end{prop}

\section{A Vanishing Theorem}\label{ch5}

In this section we start by recalling some definitions from Section 3 of
\cite{GS} and from the theory of Gr\"obner bases (see \cite{GS} and
\cite{GR}). Recall that $\L = K{\mathcal Q}/I$ where $I$ is an admissible
ideal with fixed minimal set of generators $f^2$.

\vspace{1cm}

A length-lexicographic order $>$ on the paths of ${\mathcal Q}$ is an
arbitrary linear order of both the vertices and the arrows of ${\mathcal
Q}$, so that any vertex is smaller than any path of length at least one. For
paths $p$ and $q$, both not vertices, we define $p > q$ if the length of $p$
is greater than the length of $q$. If the lengths are equal, say $p = a_1
\cdots a_t$ and $q = b_1 \cdots b_t$ where the $a_i$ and $b_i$ are arrows,
then we say $p > q$ if there is an $i, 0 \leq i \leq t - 1$, such that $a_j
= b_j$ for $j \leq i$ but $a_{i + 1} > b_{i + 1}.$

Let $f$ be an element in $K{\mathcal Q}$ written as a linear combination of
paths $\sum_{j = 1}^{s} c_j \rho_j$ with $c_j \in K \backslash \{0\}$ and
paths $\rho_j$. Following \cite{GS}, we say a path $\rho$ occurs in $f$ if
$\rho = \rho_j$ for some $j$.

Fix a length-lexicographic order on a quiver ${\mathcal Q}$. Let $f$
be a non-zero element of $K{\mathcal Q}$. Let $tip(f)$ denote the
largest path occurring in $f$. Then we define $Tip(I) = \{tip(f) | f
\in I \backslash \{0\}\}.$ Define ${\NonTip}(I)$ to be the set of
paths in $K{\mathcal Q}$ that are not in Tip(I). Note that for
vertices $v$ and $w$, $v{\NonTip}(I)w$ is a $K$-basis of paths for
$v\L w$.

\begin{defin} \cite[Definition 3.1]{GS}
The boundary of $f^2$, denoted by $Bdy(f^2)$, is defined to be the set
$$Bdy(f^2) = \{({\mathfrak o}(f^2_1),{\mathfrak t}(f^2_1)), \ldots,
({\mathfrak o}(f^2_m),{\mathfrak t}(f^2_m))\} = \{({\mathfrak
o}(x),{\mathfrak t}(x))| x \in f^2\}.$$
\end{defin}

\begin{defin} \cite[Definition 3.3]{GS}
Let ${\mathcal G}^2 = \bigcup v {\NonTip}(I)w$, where the union is taken
over all $(v, w)$ in $Bdy(f^2).$
\end{defin}
We consider now elements of ${\Hom}(Q^2, \L)$.

\begin{defin} \cite[Definition 3.4]{GS}
For $p$ in ${\mathcal G}^2$ and $x \in f^2$ with ${\mathfrak o}(x) =
{\mathfrak o}(p)$ and ${\mathfrak t}(x) = {\mathfrak t}(p)$,
define $\phi_{p_, x}: Q^2 \rightarrow \L$ to be the $\L, \L$-bimodule
homomorphism given by
$${\mathfrak o}(f^2_i) \otimes {\mathfrak t}(f^2_i) \mapsto \left\{
\begin{array}{ll}
p & \mbox{ if } f^2_i= x,\\
0 & \mbox{ otherwise}.\\
\end{array}\right.
$$
\end{defin}

Let $d_2: {\Hom}(Q^1, \L) \rightarrow {\Hom}(Q^2, \L)$ be the map induced by
$A_2$. Each element of ${\HH}^2(\L)$ may be represented by a map in
${\Hom}(Q^2, \L)$ and so is represented by a linear combination over $K$ of
maps $\phi_{p, x}$. If every $\phi_{p, x}$ is in ${\Im}\,d_2$ then
${\Hom}(Q^2, \L) = {\Im}\,d_2$ and hence ${\HH}^2(\L) = 0.$ Our strategy in
Theorem \ref{motivated thm 1} is to show that ${\HH}^2(\L) = 0$ by showing
that every $\phi_{p,x}$ is in ${\Im}\,d_2$.

First we return to \cite{GS} and modify \cite[Definition 3.6]{GS}.

\begin{defin} \label{def3.6}
Let $X$ be a set of paths in $K{\mathcal Q}$. Define

$L_0(X) = \{p \in X | \exists \mbox{ some arrow $a$ which occurs in $p$ and
}\\\hspace*{2cm} \mbox{ which does not occur in any element of } X
\backslash \{p\} \}.$\\
For $p \in L_0(X)$, we call such an $a$ an arrow associated to $p$.

Define $L_i(X)$ for $i \in {\mathbb N}$ by
$$L_i(X) = L_0(X \backslash \bigcup_{j = 0}^{i - 1} L_j(X)).$$
\end{defin}

\begin{defin} \label{defin 3.9} \cite[Definition 3.9]{GS}
Let $X$ be a set of paths in ${\NonTip}(I)$. The arrows are said to separate
$X$ if $X = \bigcup_{i \geq 0} L_i(X)$.
\end{defin}

\vspace{1cm}

Motivated by Theorem 3.10 in \cite{GS} we give a new theorem on the
vanishing of ${\HH}^2(\L)$ which we will show applies to all algebras in
Asashiba's list when $s \geq 2$. (We will consider the case $s = 1$ later.)

\begin{thm} \label{motivated thm 1}
Let $\L = K{\mathcal Q}/I$ be a finite dimensional algebra where $I$ is an
admissible ideal with minimal generating set $f^2$. With the notation of
this section, suppose that for all $(v, w) \in Bdy(f^2)$ either $v\L w = \{0
\}$ or there is some path $p$ such that $v{\NonTip}(I)w = \{p\}$. In the
case where $v \L w \neq \{0\}$ suppose further that $vf^2w = \{p - q_1,
\ldots, p - q_t\}$ for paths $q_1, \ldots, q_t.$ Thus we may write
${\mathcal G}^2 = \{p_1, \ldots, p_r\}$, where for each $i = 1, \ldots, r$,
we have non-zero paths $q_{i1}, \ldots, q_{it_{i}}$ with ${\mathfrak o}(p_i)
f^2 {\mathfrak t}(p_i) = \{p_i - q_{i1}, \ldots, p_i - q_{it_{i}}\}$.

Let $Y = \{p_1, \ldots, p_r, q_{ij}| 1 \leq i \leq r, 1 \leq j \leq t_i \}$.
Suppose that $L_0(Y) = Y$. Let $a_{ij}$ be an arrow associated to $q_{ij}$
and assume that $a_{ij}$ occurs only once in the path $q_{ij}$. Then every
element of ${\Hom}(Q^2, \L)$ is a coboundary, that is, $\phi_{p,x} \in
{\Im}\,d_2$ for all $p \in {\mathcal G}^2$ and $x \in f^2$, and thus
${\HH}^2(\L) = 0$.
\end{thm}

\begin{proof}
It is enough to show that each element $\phi_{p, x}$ of ${\Hom}(Q^2, \L)$,
where $p$ is a path in ${\mathcal G}^2$
and $x \in f^2$ with ${\mathfrak o} (x) = {\mathfrak o}(p)$ and ${\mathfrak
t} (x) = {\mathfrak t}(p)$, is a coboundary.
By hypothesis ${\mathcal G}^2 = \{p_1, \ldots, p_r\}$. Note that the paths
$p_1, \ldots, p_r$ are distinct. Consider the path $p_i$ where $i \in \{1,
\ldots, r\}$. Then by hypothesis there are vertices $v_i, w_i$ with
$v_i{\NonTip}(I)w_i = \{p_i\}$ and  $v_if^2w_i = \{p_i - q_{i1}, \ldots, p_i
- q_{it_i}\}$. Thus if $x \in f^2$ and  ${\mathfrak o} (x) = {\mathfrak
o}(p_i)$ and ${\mathfrak t} (x) = {\mathfrak t}(p_i)$ then $x \in
v_if^2w_i$. Thus $x \in \{p_i - q_{i1}, \ldots, p_i - q_{it_i}\}$. Consider
$x = p_i - q_{ij}$ where $j \in \{1, \ldots, t_i\}.$

The map $\phi_{p_i, x}: Q^2 \rightarrow \L$ is given by
$${\mathfrak o}(f^2_k) \otimes {\mathfrak t}(f^2_k) \mapsto \left\{
\begin{array}{ll}
p_i & \mbox{ if } f^2_k = x,\\
0 & \mbox{ otherwise}.\\
\end{array}\right.
$$

We have $Y = \{p_1, \ldots, p_r, q_{ij}| 1 \leq i \leq r, 1 \leq j \leq t_i
\}$ and $Y = L_0(Y)$ so $q_{ij} \in L_0(Y)$. Therefore there exists some
arrow $a_{ij}$ which occurs in $q_{ij}$ and does not occur in any element of
$Y \backslash \{q_{ij}\}$.

Define $\psi: Q^1 \rightarrow \L$ by
$${\mathfrak o}(\alpha) \otimes {\mathfrak t}(\alpha) \mapsto \left\{
\begin{array}{ll}
-a_{ij} & \mbox{ if } \alpha = a_{ij},\\
0 & \mbox{ otherwise}.\\
\end{array}\right.
$$

Now we want to show that $\psi A_2 = \phi_{p_i, x}$. Take ${\mathfrak
o}(f^2_k)\otimes {\mathfrak t}(f^2_k) \in Q^2$. We start by finding $\psi
A_2({\mathfrak o}(f^2_k) \otimes {\mathfrak t}(f^2_k))$ by considering two
cases.

\vspace{.5cm}
\noindent {\bf Case $f^2_k = x.$}

Here, we have $\psi A_2({\mathfrak o}(f^2_k) \otimes {\mathfrak t}(f^2_k)) =
  \psi A_2({\mathfrak o}(x) \otimes {\mathfrak t}(x)),$ where $x = p_i -
q_{ij}$ and $q_{ij} = \rho_1 a_{ij} \rho_2$ for paths $\rho_1, \rho_2$ such
that $a_{ij}$ does not occur in $\rho_1$ or $\rho_2$ since $a_{ij}$ occurs
only once in $q_{ij}$ by hypothesis. Let $p_i = \sigma_1 \cdots \sigma_l,
\rho_1 = \epsilon_1 \cdots \epsilon_n, \rho_2 = b_1 \cdots b_m,$ where the
$\sigma$'s, $\epsilon$'s, $b$'s are arrows.
Then $\psi A_2({\mathfrak o}(x) \otimes {\mathfrak t}(x))
=$$\psi[({\mathfrak o}(x) \otimes_{\sigma_1}(\sigma_2 \cdots \sigma_l) +
\sigma_1 \otimes_{\sigma_2}(\sigma_3 \cdots \sigma_l) + \ldots + (\sigma_1
\sigma_2 \cdots \sigma_{l - 1} \otimes_{\sigma_l} {\mathfrak t}(x) ) -
$$({\mathfrak o}(x) \otimes_{\epsilon_1}(\epsilon_2 \cdots \epsilon_n)
a_{ij} \rho_2 + \epsilon_1 \otimes_{\epsilon_2}(\epsilon_3 \cdots
\epsilon_n) a_{ij} \rho_2 + \ldots + (\epsilon_1 \epsilon_2 \cdots
\epsilon_{n-1}) \otimes_{\epsilon_n} a_{ij}\rho_2 + \rho_1 \otimes_{a_{ij}}
\rho_2 + \rho_1 a_{ij} \otimes_{b_1}(b_2 \cdots b_m) + \rho_1 a_{ij} b_1
\otimes_{b_2}(b_3 \cdots b_m) + \ldots + \rho_1 a_{ij}(b_1 b_2 \cdots b_{m -
1}) \otimes_{b_m} {\mathfrak t}(x))]$.

As $q_{ij}, p_i \in Y = L_0(Y)$ and $a_{ij}$ occurs in $q_{ij}$, we have
that $a_{ij}$ does not occur in $p_i$. So $a_{ij}$ is not equal to any of
the $\sigma$'s, $\epsilon$'s or $b$'s. Therefore

$\psi A_2({\mathfrak o}(x) \otimes {\mathfrak t}(x))$

$= -\psi(\rho_1 \otimes_{a_{ij}} \rho_2)$

$= -\rho_1 \psi({\mathfrak t}(\rho_1) \otimes_{a_{ij}} {\mathfrak
o}(\rho_2)) \rho_2$

$= -\rho_1 \psi({\mathfrak o}(a_{ij}) \otimes_{a_{ij}} {\mathfrak
t}(a_{ij})) \rho_2 $

$ = \rho_1 a_{ij} \rho_2 = q_{ij}$.

\vspace{.5cm}

\noindent {\bf Case $f^2_k \neq x.$}

We consider separately the cases ${\mathfrak o}(f^2_k) \L {\mathfrak
t}(f^2_k) = 0$ and ${\mathfrak o}(f^2_k) \L {\mathfrak t}(f^2_k) \neq 0.$

{\bf a)} If ${\mathfrak o}(f^2_k) \L {\mathfrak t}(f^2_k) = 0$ then $\psi
A_2({\mathfrak o}(f^2_k) \otimes {\mathfrak t}(f^2_k)) = {\mathfrak
o}(f^2_k)\psi A_2 ({\mathfrak o}(f^2_k) \otimes {\mathfrak
t}(f^2_k)){\mathfrak t}(f^2_k)$ \linebreak $ = 0$ as $\psi A_2({\mathfrak
o}(f^2_k) \otimes {\mathfrak t}(f^2_k)) \in \L$ and ${\mathfrak o}(f^2_k) \L
{\mathfrak t}(f^2_k) = 0.$

{\bf b)} If ${\mathfrak o}(f^2_k) \L {\mathfrak t}(f^2_k) \neq 0$ then
${\mathfrak o}(f^2_k) \L {\mathfrak t}(f^2_k) = Sp\{p_u\}$, the vector space
spanned by $p_u$, for some $1 \leq u \leq r.$ Hence $f^2_k = p_u - q_{ul}$
for some $1 \leq l \leq t_u.$

We have $L_0(Y) = Y$ so $a_{ij}$ does not occur in any element of $Y
\backslash \{q_{ij}\}.$ Suppose for contradiction that $a_{ij}$ occurs in
$q_{ul}$, so that $q_{ul} = q_{ij}$ as paths in $K{\mathcal Q}$. Then
$${\mathfrak o}(f^2_k) = {\mathfrak o}(q_{ul}) = {\mathfrak o}(q_{ij}) =
{\mathfrak o}(x)$$
and
$${\mathfrak t}(f^2_k) = {\mathfrak t}(q_{ul}) = {\mathfrak t}(q_{ij}) =
{\mathfrak t}(x).$$

Therefore, ${\mathfrak o}(f^2_k) \L {\mathfrak t}(f^2_k) = {\mathfrak o}(x)
\L {\mathfrak t}(x) = Sp\{p_i\}.$ Hence, $p_u = p_i$ by the choice of
${\mathcal G}^2.$ Therefore, $f^2_k = p_u - q_{ul} = p_i - q_{ij} = x.$ This
gives a contradiction since we assumed $f^2_k \neq x.$ Hence $a_{ij}$ does
not occur in $q_{ul}$.

Now suppose for contradiction that $a_{ij}$ occurs in $p_u$ so that $p_u =
q_{ij}$ as paths in $K{\mathcal Q}$. Then
$${\mathfrak o}(f^2_k) = {\mathfrak o}(p_u) = {\mathfrak o}(q_{ij}) =
{\mathfrak o}(x)$$
and
$${\mathfrak t}(f^2_k) = {\mathfrak t}(p_u) = {\mathfrak t}(q_{ij}) =
{\mathfrak t}(x).$$

Therefore, $Sp\{p_u\} = {\mathfrak o}(f^2_k) \L {\mathfrak t}(f^2_k) =
{\mathfrak o}(x) \L  {\mathfrak t}(x) = Sp\{p_i\}.$ Therefore, $p_u = p_i$
by the choice of ${\mathcal G}^2.$ Hence $p_i = p_u = q_{ij}$ in $K{\mathcal
Q}$. So $p_i - q_{ij} = 0$ in $K{\mathcal Q}$. This contradicts $p_i -
q_{ij}$ being a minimal generator of $I$. Therefore, $a_{ij}$ does not occur
in $p_u$.

Thus $a_{ij}$ does not occur in $f^2_k.$
So $\psi A_2({\mathfrak o}(f^2_k) \otimes {\mathfrak t}(f^2_k)) = 0.$

Hence $\psi A_2$ is the map
$${\mathfrak o}(f^2_k) \otimes {\mathfrak t}(f^2_k) \mapsto \left\{
\begin{array}{ll}
q_{ij} & \mbox{ if } f^2_k = x,\\
0 & \mbox{ otherwise}.\\
\end{array}\right.
$$
As $p_i - q_{ij} \in f^2,$ we know that $p_i = q_{ij}$ in $\L$. Hence $\psi
A_2 = \phi_{p_i, x}$. Thus $\phi_{p_i, x}$, and hence each element of
${\Hom}(Q^2, \L)$, is a coboundary. Hence ${\HH}^2(\L) = 0.$
\end{proof}

\section{Application to Standard Algebras} \label{sec4}
We now want to apply Theorem \ref{motivated thm 1} to our derived
equivalence representatives. We start by considering the standard derived
equivalence representatives, and we need minimal relations for each such
algebra in Asashiba's list.

We start with the algebra $\L = \L(D_n, s, 1)$. Note that $R(D_n, s, 1)$ for
$s \geq 1$ is not minimal.

For relations of type (i), let $\beta_{0}^{[i]} \beta_{1}^{[i]} -
\gamma_{0}^{[i]} \gamma_{1}^{[i]} \in f^2$ and $\beta_{0}^{[i]}
\beta_{1}^{[i]} - \alpha_{n - 2}^{[i]} \alpha_{n - 3}^{[i]} \cdots
\alpha_{2}^{[i]}\alpha_{1}^{[i]} \in f^2.$ All relations of type (ii) are in
$f^2$. We now consider the relations of type (iii). So $(\beta_{0}^{[i]}
\beta_{1}^{[i]} - \gamma_{0}^{[i]} \gamma_{1}^{[i]}) \gamma_{0}^{[i + 1]} =
(\beta_{0}^{[i]} \beta_{1}^{[i]} \gamma_{0}^{[i + 1]} - \gamma_{0}^{[i]}
\gamma_{1}^{[i]} \gamma_{0}^{[i + 1]}) \in I$ and $\beta_{0}^{[i]}
\beta_{1}^{[i]} \gamma_{0}^{[i + 1]} \in I$. Therefore  $\gamma_{0}^{[i]}
\gamma_{1}^{[i]} \gamma_{0}^{[i + 1]} \in I$ and is not in $f^2$. Also
$\gamma_{1}^{[i - 1]}(\beta_{0}^{[i]} \beta_{1}^{[i]} - \gamma_{0}^{[i]}
\gamma_{1}^{[i]}) = (\gamma_{1}^{[i - 1]}\beta_{0}^{[i]} \beta_{1}^{[i]} -
\gamma_{1}^{[i - 1]} \gamma_{0}^{[i]} \gamma_{1}^{[i]}) \in I$ and
$\gamma_{1}^{[i - 1]}\beta_{0}^{[i]} \beta_{1}^{[i]} \in I$. So
$\gamma_{1}^{[i - 1]} \gamma_{0}^{[i]} \gamma_{1}^{[i]} \in I$ and is not in
$f^2$.  Similarly we can show that neither $\beta_{0}^{[i]} \beta_{1}^{[i]}
\beta_{0}^{[i +1]}$ nor $\beta_{1}^{[i]} \beta_{0}^{[i + 1]} \beta_{1}^{[i
+1]}$ are in $f^2$.

Now consider ``$\alpha$-paths''. We have $\beta_{0}^{[i]} \beta_{1}^{[i]} -
\alpha_{n - 2}^{[i]} \alpha_{n - 3}^{[i]} \cdots
\alpha_{2}^{[i]}\alpha_{1}^{[i]} \in f^2.$ So $(\beta_{0}^{[i]}
\beta_{1}^{[i]} - \alpha_{n - 2}^{[i]} \alpha_{n - 3}^{[i]} \cdots
\alpha_{2}^{[i]}\alpha_{1}^{[i]}) \alpha_{n - 2}^{[i + 1]} \in I$ and
$\beta_{0}^{[i]} \beta_{1}^{[i]} \alpha_{n - 2}^{[i + 1]} \in I.$ Therefore
it follows that $\alpha_{n - 2}^{[i]} \alpha_{n - 3}^{[i]} \cdots
\alpha_{2}^{[i]}\alpha_{1}^{[i]}\alpha_{n - 2}^{[i + 1]} \in I$ and is not
in $f^2.$  Also $\alpha_{1}^{[i - 1]}(\beta_{0}^{[i]} \beta_{1}^{[i]} -
\alpha_{n - 2}^{[i]} \alpha_{n - 3}^{[i]} \cdots
\alpha_{2}^{[i]}\alpha_{1}^{[i]})$ $\in I$ and $\alpha_{1}^{[i - 1]}
\beta_{0}^{[i]} \beta_{1}^{[i]} \in I$. So $\alpha_{1}^{[i - 1]} \alpha_{n -
2}^{[i]} \alpha_{n - 3}^{[i]} \cdots \alpha_{2}^{[i]}\alpha_{1}^{[i]} \in I$
and not in $f^2$.

However, the path $\alpha_{2}^{[i]} \alpha_{1}^{[i]} \alpha_{n - 2}^{[i +
1]} \cdots \alpha_{2}^{[i + 1]}$ cannot be obtained from any other elements,
so $\alpha_{2}^{[i]} \alpha_{1}^{[i]} \alpha_{n - 2}^{[i + 1]} \cdots
\alpha_{2}^{[i + 1]} \in f^2$. In general, $\alpha_{k}^{[i]} \alpha_{k
-1}^{[i]} \cdots \alpha_{k + 1}^{[i + 1]} \alpha_{k}^{[i + 1]} \in f^2$ for
$k= \{2, \ldots, n - 3\}$. So we have the following proposition.

\vspace{.5cm}

\begin{prop}
For $\L = \L(D_n, s, 1)$ with $s \geq 1$, and
for all $i \in \{0, \ldots, s-1\}$, let
$$f_{1,1,i}^2 = \beta_{0}^{[i]} \beta_{1}^{[i]} - \gamma_{0}^{[i]}
\gamma_{1}^{[i]}, \hspace{1cm} f_{1,2,i}^2 = \beta_{0}^{[i]} \beta_{1}^{[i]}
- \alpha_{n - 2}^{[i]} \alpha_{n - 3}^{[i]} \cdots
\alpha_{2}^{[i]}\alpha_{1}^{[i]},$$       $$f_{2,1,i}^2 = \alpha_{1}^{[i]}
\beta_{0}^{[i + 1]}, \hspace{1cm} f_{2,2,i}^2 =
\alpha_{1}^{[i]}\gamma_{0}^{[i + 1]},$$
$$f_{2,3,i}^2 = \beta_{1}^{[i]}\alpha_{n - 2}^{[i + 1]},  \hspace{1cm}
f_{2,4,i}^2 = \gamma_{1}^{[i]}\alpha_{n - 2}^{[i + 1]},$$
$$f_{2,5,i}^2 = \beta_{1}^{[i]}\gamma_{0}^{[i + 1]},  \hspace{1cm}
f_{2,6,i}^2 = \gamma_{1}^{[i]}\beta_{0}^{[i + 1]} \mbox{ and }$$
$$f^2_{3,k,i} =  \alpha_k^{[i]} \cdots \alpha_1^{[i]}\alpha_{n - 2}^{[i +
1]} \cdots \alpha_k^{[i + 1]}, \mbox { for } k = \{2, \ldots, n - 3\}.$$
Then $f^2 = \{f_{1,1,i}^2, f_{1,2,i}^2, f_{2,1,i}^2, f_{2,2,i}^2,
f_{2,3,i}^2, f_{2,4,i}^2, f_{2,5,i}^2, f_{2,6,i}^2, f_{3,k,i }^2\}$ for $i =
0, \ldots, s -1$  and  $k = 2, \ldots, n - 3$ is a minimal set of relations.
\end{prop}

\vspace{.5cm}

For the rest of the algebras, we can find a minimal set of relations in a
similar way. They are given in the following propositions.

\begin{prop}\label{D_n,2}
For $\L = \L(D_n, s, 2)$ with $s \geq 2$, let,

for all $i \in \{0, \ldots, s-1\}$,
$$f_{1,1,i}^2 = \beta_{0}^{[i]} \beta_{1}^{[i]} - \gamma_{0}^{[i]}
\gamma_{1}^{[i]}, \hspace{1cm} f_{1,2,i}^2 = \beta_{0}^{[i]} \beta_{1}^{[i]}
- \alpha_{n - 2}^{[i]} \alpha_{n - 3}^{[i]} \cdots
\alpha_{2}^{[i]}\alpha_{1}^{[i]},$$       $$f_{2,1,i}^2 = \alpha_{1}^{[i]}
\beta_{0}^{[i + 1]}, \hspace{1cm} f_{2,2,i}^2 =
\alpha_{1}^{[i]}\gamma_{0}^{[i + 1]},$$
$$f_{2,3,i}^2 = \beta_{1}^{[i]}\alpha_{n - 2}^{[i + 1]},  \hspace{1cm}
f_{2,4,i}^2 = \gamma_{1}^{[i]}\alpha_{n - 2}^{[i + 1]},$$

for all $i \in \{0, \ldots, s - 2\}$,
$$f_{2,5,i}^2 = \beta_{1}^{[i]}\gamma_{0}^{[i + 1]},  \hspace{1cm}
f_{2,6,i}^2 = \gamma_{1}^{[i]}\beta_{0}^{[i + 1]},$$
$$f_{2,7,s - 1}^2 = \beta_1^{[s - 1]} \beta_0^{[0]},  \hspace{1cm}  f_{2,8,s
- 1}^2 = \gamma_1^{[s - 1]} \gamma_0^{[0]},$$

for $i \in \{0, \ldots, s-1\}$,
$$ f_{3,k,i}^2 =  \alpha_k^{[i]} \cdots \alpha_1^{[i]}\alpha_{n - 2}^{[i +
1]} \cdots \alpha_k^{[i + 1]}, \mbox { for } k = \{2, \ldots, n - 3\}.$$
Then $f^2 = \{f_{1,1,i}^2, f_{1,2,i}^2, f_{2,1,i}^2, f_{2,2,i}^2,
f_{2,3,i}^2, f_{2,4,i}^2 \mbox{ for } i = 0, \ldots, s -1 \} \cup
\{f_{2,5,i}^2,$ \linebreak $f_{2,6,i}^2 \mbox { for } i = 0, \ldots, s - 2
\} \cup \{f_{2,7,s - 1}^2, f_{2,8,s - 1}^2 \} \cup  \{f_{3,k,i }^2 \mbox{
for } i = 0, \ldots, s - 1$ and $k = 2, \ldots, n - 3\}$ is a minimal set of
relations.
\end{prop}

\vspace{.5cm}

Note that Proposition \ref{D_n,2} is for $s \geq 2$. For $s = 1$ the minimal
relations are different and are given in the next proposition.

\begin{prop}\label{D_n,1,2}
For $\L = \L(D_n, 1, 2)$, let
$$f_{1,1}^2 = \beta_{0}\beta_{1} - \gamma_{0}\gamma_{1}, \hspace{1cm}
f_{1,2}^2 = \beta_{0} \beta_{1}- \alpha_{n - 2} \alpha_{n - 3}\cdots
\alpha_{2}\alpha_{1},$$
$$f_{2,1}^2 = \alpha_{1} \beta_{0}, \hspace{1cm} f_{2,2}^2 =
\alpha_{1}\gamma_{0},$$
$$f_{2,3}^2 = \beta_{1}\alpha_{n - 2},  \hspace{1cm}  f_{2,4}^2 =
\gamma_{1}\alpha_{n - 2},$$
$$f_{2,5}^2 = \beta_1 \beta_0,  \hspace{1cm}  f_{2,6}^2 = \gamma_1 \gamma_0
\mbox{ and }$$
$$f_{3,k}^2 =  \alpha_k \cdots \alpha_1\alpha_{n - 2} \cdots \alpha_k, \mbox
{ for } k \in \{2, \ldots, n - 3\}.$$
Then $f^2 = \{f_{1,1}^2, f_{1,2}^2, f_{2,1}^2, f_{2,2}^2, f_{2,3}^2,
f_{2,4}^2, f_{2,5}^2, f_{2,6}^2, f_{3,k}^2 \mbox{ for } k = 2, \ldots, n -
3\}$ is a minimal set of relations.
\end{prop}

\vspace{.5cm}

Again for $\L(D_4,s,3)$ we separate the cases $s \geq 2$ and $s = 1$.

\begin{prop}\label{D_4,3}
For $\L = \L(D_4, s, 3)$ with $s \geq 2$, let,

for all $i \in \{0, \ldots, s - 1\}$:
$$f_{1,1,i}^2 = \beta_{0}^{[i]} \beta_{1}^{[i]} - \gamma_{0}^{[i]}
\gamma_{1}^{[i]}, \hspace{1cm} f_{1,2,i}^2 = \beta_{0}^{[i]} \beta_{1}^{[i]}
- \alpha_0^{[i]} \alpha_1^{[i]},$$
$$f_{2,1,i}^2 = \beta_1^{[i]} \alpha_0^{[i + 1]}, \hspace{1cm} f_{2,2,i}^2 =
\alpha_1^{[i]}\gamma_{0}^{[i + 1]},$$
$$f_{2,3,i}^2 = \gamma_1^{[i]}\beta_0^{[i + 1]},$$

for all $i \in \{0, \ldots, s - 2\}$:
$$f_{2,4,i}^2 = \alpha_{1}^{[i]}\beta_{0}^{[i + 1]},  \hspace{1cm}
f_{2,5,i}^2 = \beta_{1}^{[i]}\gamma_{0}^{[i + 1]},$$
$$f_{2,6,i}^2 = \gamma_1^{[i]} \alpha_0^{[i+1]},$$
$$f_{2,7,s - 1}^2 = \gamma_1^{[s - 1]} \gamma_0^{[0]}, \hspace{1cm} f_{2,8,s
- 1}^2 = \beta_1^{[s - 1]} \beta_0^{[0]},$$
$$f_{2,9,s - 1}^2 = \alpha_1^{[s - 1]} \alpha_0^{[0]};$$
$$f_{3,1,s-1}^2 = \beta_1^{[s-1]} \beta_0^{[0]} \beta_1^{[0]}, \hspace{1cm}
f_{3,2,s-1}^2 = \alpha_0^{[s-1]} \alpha_1^{[s-1]} \beta_0^{[0]},$$
$$f_{3,4,s-1}^2 = \beta_0^{[s-1]} \beta_1^{[s-1]} \gamma_0^{[0]},
\hspace{1cm} f_{3,5,s-1}^2 = \alpha_1^{[s-1]} \beta_0^{[0]} \beta_1^{[0]}
\mbox{ and }$$
$$f_{3,6,s-1}^2 = \beta_1^{[s-1]} \gamma_0^{[0]} \gamma_1^{[0]},
\hspace{1cm} f_{3,7,s-1}^2 = \gamma_1^{[s-1]} \alpha_0^{[0]}
\alpha_1^{[0]}.$$
Then $f^2 = \{f_{1,1,i}^2, f_{1,2,i}^2, f_{2,1,i}^2, f_{2,2,i}^2,
f_{2,3,i}^2, \mbox{ for } i = 0, \ldots, s -1 \} \cup \{f_{2,4,i}^2,
f_{2,5,i}^2,$ \linebreak $f_{2,6,i}^2 \mbox { for } i = 0, \ldots, s - 2 \}
\cup \{f_{2,7,s - 1}^2, f_{2,8,s - 1}^2, f_{2,9,s - 1}, f_{3,1,s-1}^2,
f_{3,2,s-1}^2, f_{3,3,s-1}^2,$ \linebreak $f_{3,4,s-1}^2, f_{3,5,s-1}^2,
f_{3,6,s-1}^2\}$ is a minimal set of relations.
\end{prop}

\vspace{.5cm}

\begin{prop}\label{D_4,1,3}
For $\L = \L(D_4, 1, 3)$, let
$$f_{1,1}^2 = \beta_0\beta_1 - \gamma_0\gamma_1, \hspace{1cm} f_{1,2}^2 =
\beta_0 \beta_1 - \alpha_0 \alpha_1,$$
$$f_{2,1}^2 = \beta_1\alpha_0, \hspace{1cm} f_{2,2}^2 = \alpha_1\gamma_0,$$
  $$f_{2,3}^2 = \gamma_1\beta_0,$$
$$f_{2,4}^2 = \gamma_1 \gamma_0, \hspace{1cm} f_{2,5}^2 = \beta_1 \beta_0
\mbox{ and }$$
$$f_{2,6}^2 =  \alpha_1\alpha_0.$$
Then $f^2 = \{f_{1,1}^2, f_{1,2}^2, f_{2,1}^2, f_{2,2}^2, f_{2,3}^2,
f_{2,4}^2, f_{2,5}^2, f_{2,6}^2\}$ is a minimal set of relations.
\end{prop}

\vspace{1cm}

\begin{prop} \label{D_{3m},1}
For the standard algebra $\L = \L(D_{3m}, s/3, 1)$ with $s \geq 1$, for all
$i \in \{1, \ldots, s\}$, let
$$f_{1,i}^2 = \beta_i\beta_{i+1} - \alpha_1^{[i]} \cdots \alpha_m^{[i]},
\hspace{1cm} f_{2,i}^2 = \alpha_m^{[i]} \alpha_1^{[i+2]},$$
$$f_{3,i,j}^2 = \alpha_j^{[i]} \cdots \alpha_m^{[i]} \beta_{i+2}
\alpha_1^{[i+3]} \cdots \alpha_j^{[i+3]} \mbox{ for all } j \in \{2, \ldots,
m - 1\}.$$
Then $f^2 = \{f_{1,i}^2, f_{2,i}^2, f_{3,i,j}^2 \mbox{ for } j = 2, \ldots,
m-1 \mbox{ and } i = 1, \ldots, s \}$ is a minimal set of relations.
\end{prop}

\vspace{1cm}

\begin{prop} \label{E_n,1}
For $\L = \L(E_n, s, 1)$ with $s \geq 1$ and for all $i \in \{0, \ldots, s -
1\}$, let
$$f_{1,1,i}^2 = \beta_3^{[i]} \beta_2^{[i]} \beta_1^{[i]} - \gamma_2^{[i]}
\gamma_{1}^{[i]}, \hspace{1cm} f_{1,2,i}^2 = \beta_3^{[i]} \beta_2^{[i]}
\beta_1^{[i]} - \alpha_{n-3}^{[i]} \alpha_{n-4}^{[i]} \cdots \alpha_2^{[i]}
\alpha_1^{[i]},$$
$$f_{2,1,i}^2 = \alpha_1^{[i]} \beta_3^{[i + 1]}, \hspace{1cm} f_{2,2,i}^2 =
\alpha_1^{[i]}\gamma_2^{[i + 1]},$$
$$f_{2,3,i}^2 = \beta_1^{[i]}\alpha_{n-3}^{[i + 1]}, \hspace{1cm}
f_{2,4,i}^2 = \beta_1^{[i]}\gamma_2^{[i + 1]},$$
$$f_{2,5,i}^2 = \gamma_1^{[i]}\alpha_{n-3}^{[i + 1]}, \hspace{1cm}
f_{2,6,i}^2 = \gamma_1^{[i]} \beta_3^{[i+1]},$$
$$f_{3,k,i}^2 = \alpha_k^{[i]} \alpha_{k-1}^{[i]} \cdots
\alpha_{k+1}^{[i+1]}\alpha_k^{[i+1]} \mbox{ for } k \in \{2, \ldots, n-4\}
\mbox{ and }$$
$$f_{4,i}^2 = \beta_2^{[i]} \beta_1^{[i]} \beta_3^{[i]} \beta_2^{[i+1]}.$$
Then $f^2 = \{f_{1,1,i}^2, f_{1,2,i}^2, f_{2,1,i}^2, f_{2,2,i}^2,
f_{2,3,i}^2, f_{2,4,i}^2, f_{2,5,i}^2, f_{2,6,i}^2, f_{3,k,i}^2 \mbox{ for }
k \in \{2, \ldots,$ \linebreak $n-4\}, f_{4,i}^2\}$ is a minimal set of
relations.
\end{prop}

\vspace{1cm}

Finally, for the algebras of type $E_6$ we have 2 cases to consider.

\begin{prop}\label{E_6,2}
For $\L = \L(E_6, s, 2)$ with $s \geq 2$, let,\\
for all $i \in \{0, \ldots, s - 1\}$:
$$f_{1,1,i}^2 = \beta_3^{[i]} \beta_2^{[i]} \beta_1^{[i]} - \gamma_2^{[i]}
\gamma_{1}^{[i]}, \hspace{1cm} f_{1,2,i}^2 = \beta_3^{[i]} \beta_2^{[i]}
\beta_1^{[i]} - \alpha_3^{[i]} \alpha_2^{[i]} \alpha_1^{[i]},$$
$$f_{2,1,i}^2 = \gamma_1^{[i]} \alpha_3^{[i + 1]}, \hspace{1cm} f_{2,2,i}^2
= \alpha_1^{[i]}\beta_3^{[i + 1]},$$
$$f_{2,3,i}^2 = \alpha_1^{[i]}\gamma_2^{[i + 1]}, \hspace{1cm} f_{2,4,i}^2 =
\beta_1^{[i]}\gamma_2^{[i + 1]},$$

and for all $i \in  \{0, \ldots, s-2\}$:
$$f_{2,5,i}^2 = \alpha_1^{[i]}\beta_3^{[i + 1]}, \hspace{1cm} f_{2,6,i}^2 =
\beta_1^{[i]} \alpha_3^{[i+1]},$$

$$f_{2,7,s-1}^2 = \alpha_1^{[s-1]} \alpha_3^{[0]}, \hspace{1cm}
f_{2,8,s-1}^2 = \beta_1^{[s-1]} \beta_3^{[0]}$$

$$f_{3,1,i}^2 = \alpha_2^{[i]} \alpha_1^{[i]} \alpha_3^{[i+1]}
\alpha_2^{[i+1]}, \hspace{1cm} f_{3,2,i}^2 = \beta_2^{[i]} \beta_1^{[i]}
\beta_3^{[i+1]} \beta_2^{[i+1]},$$

$$f_{3,3,s-1}^2 = \alpha_2^{[s-1]} \alpha_1^{[s-1]} \beta_3^{[0]}
\beta_2^{[0]}, \hspace{1cm} f_{3,4,s-1}^2 = \beta_2^{[s-1]} \beta_1^{[s-1]}
\alpha_3^{[0]} \alpha_2^{[0]}.$$
Then $f^2 = \{f_{1,1,i}^2, f_{1,2,i}^2, f_{2,1,i}^2, f_{2,2,i}^2,
f_{2,3,i}^2, f_{2,4,i}^2, \mbox{ for } i = 0, \ldots s-1\} \cup
\{f_{2,5,i}^2, f_{2,6,i}^2$, $\mbox{ for } i = 0, \ldots, s-2 \} \cup
\{f_{2,7,s-1}^2, f_{2,8,s-1}^2\} \cup \{f_{3,1,i}^2, f_{3,2,i}^2, \mbox{ for
} i = 0, \ldots, s-2\} \cup \{f_{3,3,s-1}^2, f_{3,4,s-1}\}$ is a minimal set
of relations.
\end{prop}

\vspace{1cm}

\begin{prop}\label{E_6,1,2}
For $\L = \L(E_6, 1, 2)$, let
$$f_{1,1}^2 = \beta_3 \beta_2 \beta_1 - \gamma_2 \gamma_1, \hspace{1cm}
f_{1,2}^2 = \beta_3 \beta_2 \beta_1 - \alpha_3 \alpha_2 \alpha_1,$$
$$f_{2,1}^2 = \gamma_1 \alpha_3, \hspace{1cm} f_{2,2}^2 = \gamma_1
\beta_3,$$
$$f_{2,3}^2 = \alpha_1 \gamma_2, \hspace{1cm} f_{2,4}^2 = \beta_1
\gamma_2,$$
$$f_{2,5}^2 = \alpha_1 \alpha_3, \hspace{1cm} f_{2,6}^2 = \beta_1 \beta_3,$$
   $$f_{3,1}^2 = \alpha_2 \alpha_1 \beta_3 \beta_2, \hspace{1cm} f_{3,2}^2 =
\beta_2 \beta_1 \alpha_3 \alpha_2.$$
        Then $f^2 = \{f_{1,1}^2, f_{1,2}^2, f_{2,1}^2, f_{2,2}^2, f_{2,3}^2,
f_{2,4}^2, f_{2,5}^2, f_{2,6}^2, f_{3,1}^2, f_{3,2}^2\}$ is a minimal set of
relations.
\end{prop}

\vspace{1cm}

We now apply Theorem \ref{motivated thm 1} to the self-injective algebras of
type $D_n$ and $E_{6,7,8}$ using Propositions \ref{D_n,2}, \ref{D_n,1,2},
\ref{D_4,3}, \ref{D_4,1,3}, \ref{D_{3m},1}, \ref{E_n,1}, \ref{E_6,2} and
\ref{E_6,1,2}.

For example consider the algebra $\L(D_n, s, 2)$ for $s \geq 2$. Fix an
order on the vertices and the arrows:\\
$\alpha_{n-2}^{[0]} > \alpha_{n-3}^{[0]} > \cdots > \alpha_1^{[0]} >
\gamma_0^{[0]} > \gamma_1^{[0]} > \beta_0^{[0]} > \beta_1^{[0]} >$
$\alpha_{n-2}^{[1]} > \cdots > \beta_1^{[1]} >$ $ \cdots >
\alpha_{n-2}^{[s-1]} > \cdots > \beta_1^{[s-1]}$\\
and $\beta_1^{[s-1]} > e_{1,0} > e_{n-2,0} > \cdots > e_{1,1} > e_{n,0} >
e_{n-1,0} > \cdots > e_{1,s-1} > e_{n-2,s-1} > \cdots > e_{n,s-1} >
e_{n-1,s-1}.$

Then $tip(f^2_{1,1,i}) = tip(\beta_0^{[i]} \beta_1^{[i]} - \gamma_0^{[i]}
\gamma_1^{[i]}) = \gamma_0^{[i]} \gamma_1^{[i]}$ and $tip(f^2_{1,2,i}) =
tip(\beta_0^{[i]} \beta_1^{[i]} - \alpha_{n - 2}^{[i]} \alpha_{n - 3}^{[i]}
\cdots \alpha_2^{[i]}\alpha_1^{[i]}) = \alpha_{n - 2}^{[i]} \alpha_{n -
3}^{[i]} \cdots \alpha_2^{[i]}\alpha_1^{[i]}$ for $i = 0, \ldots, s-1.$ For
all other $f^2_j \in f^2$ with $f^2_j \neq f^2_{1,1,i}, f^2_{1,2,i}$ we know
that $f^2_j$ is a path in $K{\mathcal Q}$ so $tip(f^2_j) = f^2_j.$ In these
cases $\mo(f^2_j){\NonTip}(I)\mt(f^2_j) = \{0\}$. Let $v_i =
\mo(f^2_{1,1,i}) = \mo(f^2_{1,2,i})$ and let $w_i = \mt(f^2_{1,1,i}) =
\mt(f^2_{1,2,i})$ for $i = 0, \ldots, s-1.$ Then $(v_i, w_i) \in Bdy(f^2)$
and $v_i{\NonTip}(I)w_i = \{\beta_0^{[i]} \beta_1^{[i]}\}$ for all $i = 0,
\ldots, s-1.$ So let $p^{[i]} = \beta_0^{[i]} \beta_1^{[i]}$ for $i = 0,
\ldots, s-1.$ Then $v_if^2w_i = \{\beta_0^{[i]} \beta_1^{[i]} -
\gamma_0^{[i]} \gamma_1^{[i]}, \beta_0^{[i]} \beta_1^{[i]} - \alpha_{n -
2}^{[i]} \alpha_{n - 3}^{[i]} \cdots \alpha_2^{[i]}\alpha_1^{[i]}\} =
\{p^{[i]} - q^{[i]}_1, p^{[i]} - q^{[i]}_2\}$, where $q_1^{[i]} =
\gamma_0^{[i]} \gamma_1^{[i]}$, $q_2^{[i]} = \alpha_{n - 2}^{[i]} \alpha_{n
- 3}^{[i]} \cdots \alpha_2^{[i]}\alpha_1^{[i]}.$ With the notation of
Theorem \ref{motivated thm 1}, ${\mathcal G}^2 = \{\beta_0^{[i]}
\beta_1^{[i]} \ |\ i = 0, \ldots, s-1\}$ and $Y = \{\beta_0^{[i]}
\beta_1^{[i]}, \gamma_0^{[i]} \gamma_1^{[i]}, \alpha_{n - 2}^{[i]} \alpha_{n
- 3}^{[i]} \cdots \alpha_2^{[i]}\alpha_1^{[i]}$ $\|\ i = 0, \ldots, s-1\} =
L_0(Y).$ Choose $a_1^{[i]} = \gamma_0^{[i]}$ and $a_2^{[i]} = \alpha_{n -
2}^{[i]}$ so that $a_1^{[i]}$ and $a_2^{[i]}$ are arrows associated to
$q_1^{[i]}$ and $q_2^{[i]}$ respectively, and $a_j^{[i]}$ occurs once in
$q_j^{[i]}$ for $j = 1, 2$. Then by applying Theorem \ref{motivated thm 1},
every element of ${\Hom}(Q^2, \L)$ is a coboundary and so ${\HH}^2(\L) = 0.$

Similar arguments give the following corollary.

\begin{cor} \label{cor 1}
Suppose $s \geq 2$. Let $\L$ be one of the standard algebras $\L(D_n, s,
1),$ $\L(D_n, s, 2)$ for $n \geq 4$, $\L(D_4, s, 3)$, $\L(D_{3m}, s/3, 1)$
with $m \geq 2, 3 \nmid s,$ $\L(E_n, s, 1)$ with $n \in \{6, 7, 8\}$ or
$\L(E_6, s, 2)$. Then ${\HH^2}(\L) = 0$.
\end{cor}

\begin{remark}
Theorem \ref{motivated thm 1} does not apply if $s = 1$ since in this case
there is some $(v, w) \in Bdy(f^2)$ with $\dim\, v \L w > 1$.
\end{remark}

\section{${\HH}^2(\L)$ for the standard self-injective algebras of finite
representation type}\label{sec5}

In this section we determine ${\HH}^2(\L)$ for the standard algebras
$\L(D_n, s, 1)$, $\L(D_n, s, 2),$ $\L(D_4, s, 3),$ $\L(D_{3m}, s/3, 1),$
$\L(E_n, s, 1),$ $\L(E_6, s, 1)$ when $s = 1.$ A sketch of the proof is
given in each type. We start with $\L(D_n, s, 2)$ since ${\HH}^2(\L) \neq 0$
in this case.

\begin{thm}\label{sumch6}
For $\L = \L(D_n, 1, 2)$ we have $\dim\,{\HH}^2(\L) = 1.$
\end{thm}

\begin{proof}
For $\L = \L(D_n, 1, 2)$ we label the quiver $Q(D_n, 1)$ as follows:
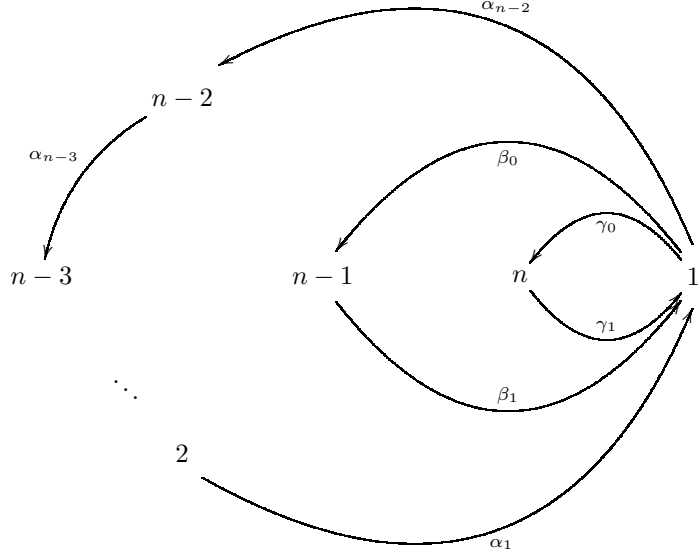
\vspace{1cm}
$$\xymatrix{
& & & & & & & & \\
& & {n-2}\ar@/_1pc/[ddl]_{\alpha_{n-3}} & & & & & & \\
& & & & & & & & \\
& {n-3} & & {n-1}\ar@<-.7ex>@/_4pc/[rrrr]^{\beta_1} & &
n\ar@/_2pc/[rr]^{\gamma_1} & &
1\ar@/_2pc/[ll]^{\gamma_0}\ar@<-.7ex>@/_4pc/[llll]^{\beta_0}\ar@<-1.5ex>@/_5pc/[uulllll]_{\alpha_{n-2}}
& \\
& \ar@{}[dr]^\ddots & & & & & & & & \\
& & 2 \ar@<-1.5ex>@/_5pc/[uurrrrr]_{\alpha_1} & & & & & &\\
}$$
\vspace{1cm}

The set $f^2$ of minimal relations was given in Proposition \ref{D_n,1,2}.
Recall that the projective $Q^3 = \bigoplus_{y \in f^3} \L \mo(y) \otimes
\mt(y) \L$ $= (\L e_1 \otimes e_{n-3}\L) \oplus (\L  e_1 \otimes e_{n -
2}\L) \oplus (\L e_1 \otimes e_{n - 1}\L) \oplus (\L  e_1 \otimes e_n \L)
\oplus (\L e_2 \otimes e_1\L) \oplus (\L e_{n-1} \otimes e_1 \L) \oplus (\L
e_n \otimes e_1\L) \oplus \bigoplus_{m = 3}^{n - 2} (\L e_m \otimes e_{m -
2}\L).$
(We note that the projective $Q^3$ is also described in \cite{H} although
Happel gives no description of the maps in the $\L, \L$-projective
resolution of $\L$.) Following \cite{GS}, and with the notation introduced
in Section \ref{ch2}, we may choose the set $f^3$ to consist of the
following elements:
$$\{f^3_{1,1}, f^3_{1,2}, f^3_{1,3}, f^3_{1,4}, f^3_2, f^3_{n-1}, f^3_n,
f^3_3, f^3_m\}, \mbox{ with } m \in \{4, \ldots, n - 2\} \mbox{  where }$$
$$\begin{array}{l c l c l l l}
f^3_{1,1} & = & f^2_{1,2} \alpha_{n - 2} \alpha_{n - 3} & = & \beta_0
f^2_{2,3} \alpha_{n - 3} - \alpha_{n - 2} f^2_{3,n-3} & \in & e_1 K{\mathcal
Q} e_{n-3},\\

f^3_{1,2} & = & f^2_{1,1} \alpha_{n - 2} & = & \beta_0 f^2_{2,3} - \gamma_0
f^2_{2,4} & \in & e_1 K{\mathcal Q} e_{n-2},\\

f^3_{1,3} & = & f^2_{1,2} \beta_0 & = & \beta_0 f^2_{2,5} -
\alpha_{n-2}\cdots \alpha_2 f^2_{2,1} & \in & e_1 K{\mathcal Q} e_{n-1},\\

f^3_{1,4} & = & f^2_{1,1} \gamma_0 - f^2_{1,2} \gamma_0 & = & \alpha_{n-2}
\cdots \alpha_2 f^2_{2,2} - \gamma_0 f^2_{2,6} & \in & e_1 K{\mathcal Q}
e_n,\\

f^3_2 & = & f^2_{2,1} \beta_1 - f^2_{2,2} \gamma_1 & = & \alpha_{1}
f^2_{1,1} & \in & e_2 K{\mathcal Q} e_1,\\

f^3_{n-1} & = & f^2_{2,5} \beta_1 - f^2_{2,3} \alpha_{n - 3} \cdots \alpha_1
& = & \beta_1 f^2_{1,2} & \in & e_{n-1} K{\mathcal Q} e_1,\\

f^3_n & = & f^2_{2,4} \alpha_{n -3} \cdots \alpha_1 - f^2_{2,6} \gamma_1 & =
& \gamma_1 f^2_{1,1} - \gamma_1 f^2_{1,2} & \in & e_n K{\mathcal Q} e_1,\\

f^3_3 & = & f^2_{3, 2} \alpha_1 & = & \alpha_2 f^2_{2,1}\beta_1 - \alpha_2
\alpha_1 f^2_{1,2} & \in & e_3 K{\mathcal Q} e_1,\\

f^3_m & = & f^2_{3, m -1} \alpha_{m -2} & = & \alpha_{m - 1} f^2_{3,m - 2} &
\in & e_m K{\mathcal Q} e_{m-2} \\

& & & & & & \mbox{ for } m \in \{4, \ldots, n - 2\}. \\
\end{array}$$

\vspace{1cm}

We know that ${\HH}^2(\L) = {\Ker}\,d_3 / {\Im}\,d_2$. First we will find
${\Im}\,d_2$. Let $f \in {\Hom}(Q^1, \L)$ and so write
$$f(e_1 \otimes_{\beta_0} e_{n - 1}) = c_1 \beta_0, \hspace{1cm} f(e_{n - 1}
\otimes_{\beta_1} e_1) = c_2 \beta_1,$$
$$f(e_1 \otimes_{\gamma_0} e_n) = c_3 \gamma_0, \hspace{1cm} f(e_n
\otimes_{\gamma_1} e_1) = c_4 \gamma_1,$$
$$f(e_1 \otimes_{\alpha_{n-2}} e_{n-2}) = d_{n-2} \alpha_{n-2}$$ and
$$f(e_{l + 1} \otimes_{\alpha_{l}} e_{l}) = d_{l} \alpha_{l} \mbox{ for } l
\in \{1, \ldots, n - 3\},$$
where $c_1, c_2, c_3, c_4, d_l \in K \mbox{ for } l \in \{1, \ldots, n-2\}.$

Now we find $fA_2 = d_2f$. We have

$fA_2(e_1 \otimes_{f^2_{1,1}} e_{1}) = f(e_{1} \otimes_{\beta_0} e_{n -1})
\beta_1 - f(e_{1} \otimes_{\gamma_0} e_{n})\gamma_1 + \beta_0 f(e_{n - 1}
\otimes_{\beta_1} e_{1}) - \gamma_0 f(e_{n} \otimes_{\gamma_1} e_{1}) =
c_{1} \beta_0 \beta_1 - c_{3} \gamma_0 \gamma_1 +  c_{2} \beta_0 \beta_1 -
c_{4} \gamma_0 \gamma_1 = (c_{1} - c_{3} + c_{2} - c_{4}) \beta_0 \beta_1.$

Also
$fA_2(e_{1} \otimes_{f^2_{1,2}} e_{1}) =  f(e_{1} \otimes_{\beta_0}
e_{n-1})\beta_1 + \beta_0 f(e_{n-1} \otimes_{\beta_1} e_{1}) - f(e_{1}
\otimes_{\alpha_{n- 2}} e_{n-2}) \alpha_{n-3}\cdots \alpha_1 -  \alpha_{n-2}
f(e_{n-2} \otimes_{\alpha_{n-3}} e_{n-3})\alpha_{n-4}\cdots \alpha_1 -
\ldots - \alpha_{n-2} \cdots \alpha_2 f(e_{2} \otimes_{\alpha_1} e_{1}) =
c_{1} \beta_0 \beta_1 + c_{2} \beta_0 \beta_1 - d_{n-2} \alpha_{n-2} \cdots
\alpha_1 -  \ldots - d_{1} \alpha_{n-2} \cdots \alpha_2 \alpha_1 = (c_{1} +
c_{2} - d_{n-2} -  \ldots - d_{1})  \beta_0 \beta_1.$

By direct calculation, we may show that $fA_2$ is given by
$$fA_2(e_1 \otimes_{f^2_{1,1}} e_{1}) = (c_{1} - c_{3} + c_{2} - c_{4})
\beta_0 \beta_1 = c' \beta_0 \beta_1,$$
$$fA_2(e_{1} \otimes_{f^2_{1,2}} e_{1}) = (c_{1} + c_{2} - d_{n-2} -  \ldots
- d_{1})  \beta_0 \beta_1 = c'' \beta_0 \beta_1$$
for some $c', c'' \in K$ and
$$fA_2(\mo(f^2_j) \otimes \mathfrak{t}(f^2_j)) = 0$$
for all $f^2_j \neq f^2_{1,1}, f^2_{1,2}$
. So $\dim\,{\Im}\,d_2 = 2$.

Now we determine ${\Ker}\,d_3$. Let $h \in {\Ker}\,d_3$, so  $h \in
{\Hom}(Q^2, \L)$ and $d_3h = 0$. Then $h: Q^2 \rightarrow \L$ is given by

$$h(e_1 \otimes_{f_{1,1}^2} e_1) = c_1 e_1 + c_2 \beta_0\beta_1,$$
$$h(e_1 \otimes_{f_{1,2}^2} e_1) = c_3 e_1 + c_4 \beta_0\beta_1,$$
$$h(\mo(f^2_{2,j}) \otimes_{f^2_{2,j}} \mt(f^2_{2,j})) = 0, \mbox{ for } j
\in \{1, \ldots, 4\},$$
$$h(e_{n - 1} \otimes_{f_{2,5}^2} e_{n-1}) = c_5 e_{n-1},$$
$$h(e_n \otimes_{f_{2,6}^2} e_n) = c_6 e_n \mbox{ and }$$
$$h(\mo(f^2_{3,k}) \otimes_{f^2_{3,k}} \mt(f^2_{3,k})) = d_k \alpha_k,
\mbox{ for } k \in \{2, \ldots, n - 3\}$$
for some $c_1, \ldots, c_6, d_k \in K$ for $k \in \{2, \ldots, n-3\}.$

Then $hA_3(e_1 \otimes_{f^3_{1,1}} e_{n - 3}) = h(e_1 \otimes_{f^2_{1,2}}
e_1) \alpha_{n - 2} \alpha_{n - 3} - \beta_0 h(e_{n - 1}\otimes_{f^2_{2,3}}
e_{n - 2}) \alpha_{n - 3}$ \linebreak $+ \alpha_{n-2} h(e_{n-2}
\otimes_{f^2_{3,n-3}} e_{n - 3}) = (c_3 e_1 + c_4 \beta_0\beta_1)
\alpha_{n-2} \alpha_{n-3} - 0 + d_{n-3}\alpha_{n-2}\alpha_{n-3} = (c_3 +
d_{n-3}) \alpha_{n-2}\alpha_{n-3}.$ As $h \in {\Ker}\,d_3$ we have $c_3 +
d_{n-3} = 0.$

In a similar way, by considering $hA_3 (\mo(f^3_{l}) \otimes_{f^3_{l}}
\mt(f^3_{l}))$ for all $f^3_l \in f^3$, $f^3_l \neq f^3_{1,1}$ it follows
that $h$ is given by
$$h(e_1 \otimes_{f_{1,1}^2} e_1) = c_2 \beta_0\beta_1,$$
$$h(e_1 \otimes_{f_{1,2}^2} e_1) = c_3 e_1 + c_4 \beta_0\beta_1,$$
$$h(\mo(f^2_{2,j}) \otimes_{f^2_{2,j}} \mt(f^2_{2,j})) = 0, \mbox{ for } j
\in \{1, \ldots, 4\},$$
$$h(e_{n-1} \otimes_{f^2_{2,5}} e_{n-1}) = c_3 e_{n-1}$$
$$h(e_n \otimes_{f^2_{2,6}} e_n) = c_3 e_n \mbox{ and }$$
$$h(\mo(f^2_{3,k}) \otimes_{f^2_{3,k}} \mt(f^2_{3,k})) = -c_3 \alpha_k,
\mbox{ for } k \in \{2, \ldots, n - 3\}$$
for some $c_2, c_3, c_4 \in K$. Hence $\dim\,{\Ker}\,d_3 = 3.$

Therefore  $\dim\,{\HH}^2(\L) = \dim\,{\Ker}\,d_3 - \dim\,{\Im}\,d_2 = 3 - 2
= 1.$
\end{proof}

\vspace{.5cm}

\begin{bit}\label{basisD_n}
{\it A basis for ${\HH}^2(\L)$ for $\L = \L(D_n, 1, 2)$.}
\end{bit}
Let $\eta$ be the map in ${\Ker}\,d_3$ given by
$$\begin{array}{rcl}
e_1 \otimes_{f_{1,2}^2} e_1 & \mapsto & e_1,\\
e_{n-1} \otimes_{f^2_{2,5}} e_{n-1} & \mapsto & e_{n-1},\\
e_n \otimes_{f^2_{2,6}} e_n & \mapsto & e_n,\\
\mo(f^2_{3,k}) \otimes_{f^2_{3,k}} \mt(f^2_{3,k}) & \mapsto & - \alpha_k,
\mbox{ for } k \in \{2, \ldots, n - 3\},\\
\mbox{else } & \mapsto & 0.
\end{array}$$

Clearly, $\eta$ is a non-zero map. Suppose for contradiction that $\eta \in
{\Im}\,d_2$. Then by the definition of $\eta$, we have $\eta(e_n
\otimes_{f^2_{2,6}} e_n) = e_n.$ On the other hand, $\eta(e_n
\otimes_{f^2_{2,6}} e_n) = fA_2(e_n \otimes_{f^2_{2,6}} e_n)$ for some $f
\in {\Hom}(Q^1, \L)$. So $\eta(e_n \otimes_{f^2_{2,6}} e_n) = 0.$ So we have
a contradiction. Therefore  $\eta \not \in {\Im}\,d_2$.

Thus $\eta + {\Im}\,d_2$ is a non-zero element of ${\HH}^2(\L)$ and the set
$\{\eta + {\Im}\,d_2\}$ is a basis of ${\HH}^2(\L)$.

\vspace{.5cm}

\begin{thm}\label{sumch4}
For $\L = \L(D_n, 1, 1)$ with $n \geq 4,$ we have ${\HH}^2(\L) = 0.$
\end{thm}

\begin{proof}
With the quiver ${\mathcal Q}(D_n, 1)$ as in Theorem \ref{sumch6} and direct
calculations for $s = 1$ we choose the set $f^3$ to consist of the following
elements:
$$\{f^3_{1,1}, f^3_{1,2}, f^3_{1,3}, f^3_{1,4}, f^3_2, f^3_{n-1}, f^3_n,
f^3_3, f^3_m\}, \mbox{ with } m \in \{4, \ldots, n - 2\} \mbox{  where }$$
$$\begin{array}{l c l c l}
f^3_{1,1} & = & f^2_{1,2} \alpha_{n - 2}  \alpha_{n-3}    & = &  \beta_0
f^2_{2,3} \alpha_{n-3} - \alpha_{n-2} f^2_{3,n-3} \in e_1 K{\mathcal Q}
e_{n-3},\\

f^3_{1,2} & = & f^2_{1,1} \alpha_{n - 2} & = & \beta_0 f^2_{2,3} - \gamma_0
f^2_{2,4} \in e_1 K{\mathcal Q} e_{n-2},\\

f^3_{1,3} & = & f^2_{1,1} \beta_0 - f^2_{1,2} \beta_0 & = & \alpha_{n-2}
\cdots \alpha_2 f^2_{2,1} - \gamma_0 f^2_{2,6} \in e_1 K{\mathcal Q}
e_{n-1},\\

f^3_{1,4} & = & f^2_{1,2} \gamma_0 & = & \beta_0 f^2_{2,5} - \alpha_{n-2}
\cdots \alpha_2 f^2_{2,2} \in e_1 K{\mathcal Q} e_n,\\

f^3_2 & = & f^2_{2,1} \beta_1 - f^2_{2,2} \gamma_1 & = & \alpha_{1}
f^2_{1,1} \in e_2 K{\mathcal Q} e_1,\\

f^3_{n-1} & = & f^2_{2,3} \alpha_{n-3} \cdots \alpha_1 - f^2_{2,5} \gamma_1
& = & \beta_1 f^2_{1,1} - \beta_1 f^2_{1,2} \in e_{n-1} K{\mathcal Q} e_1,\\

f^3_n & = & f^2_{2,6} \beta_1 - f^2_{2,4} \alpha_{n -3} \cdots \alpha_1 & =
& \gamma_{1} f^2_{1,2} \in e_n K{\mathcal Q} e_1,\\

f^3_3 & = & f^2_{3, 2} \alpha_1 & = & \alpha_2 f^2_{2,1}\beta_1 - \alpha_2
\alpha_1 f^2_{1,2} \in e_3 K{\mathcal Q}e_1,\\

f^3_m & = & f^2_{3, m -1} \alpha_{m -2} & = & \alpha_{m - 1} f^2_{3,m - 2}
\in e_m K{\mathcal Q} e_{m-2}\\
&&&& \mbox{ for } m \in \{4, \ldots, n - 2\}.
\end{array}$$
Then it is straightforward to show that $\dim\,{\Im}\,d_2 =
\dim\,{\Ker}\,d_3 = 2$ and so ${\HH}^2(\L) = 0.$
\end{proof}

\vspace{.5cm}

\begin{thm}\label{sumch7}
For $\L = \L(D_4, 1, 3)$ we have ${\HH}^2(\L) = 0.$
\end{thm}

\begin{proof}
We have the quiver ${\mathcal Q}(D_4, 1)$ as in Theorem \ref{sumch6} with $n
= 4$ and, following Asashiba in \cite{A}, write $\alpha_0$ for $\alpha_2$.
By direct calculation we choose  the following set $f^3 = \{f^3_{1,1},
f^3_{1,2}, f^3_{1,3}, f^3_2, f^3_3, f^3_4\}$ where
$$\begin{array}{l c l c l}
f^3_{1,1} & = & f^2_{1,1} \gamma_0 - f^2_{1,2} \gamma_0 & = &  \alpha_0
f^2_{2,2} - \gamma_0 f^2_{2,4} \in e_1 K{\mathcal Q} e_4,\\
f^3_{1,2} & = & f^2_{1,1} \beta_0 & = & \beta_0 f^2_{2,5} - \gamma_0
f^2_{2,3} \in e_1 K{\mathcal Q} e_3,\\
f^3_{1,3} & = & f^2_{1,2} \alpha_0 & = & \beta_0 f^2_{2,1} - \alpha_0
f^2_{2,6}\in e_1K{\mathcal Q} e_2,\\
f^3_2 & = & f^2_{2,6} \alpha_1 - f^2_{2,2} \gamma_1 & = & \alpha_1 f^2_{1,1}
- \alpha_1 f^2_{1,2} \in e_2 K{\mathcal Q} e_1,\\
f^3_3 & = & f^2_{2,5} \beta_1 - f^2_{2,1} \alpha_1 & = & \beta_1 f^2_{1,2}
\in e_3 K{\mathcal Q} e_1,\\
f^3_4 & = & f^2_{2,3} \beta_1 - f^2_{2,4} \gamma_1 & = & \gamma_1 f^2_{1,1}
\in e_4 K{\mathcal Q} e_1.\\
\end{array}$$
We can then show that $\dim\,{\Ker}\,d_3 = \dim\,{\Im}\,d_2 = 2$ and so
${\HH}^2(\L) = 0.$
\end{proof}

\vspace{.5cm}

\begin{thm}\label{sumch8}
For the standard algebra $\L = \L(D_{3m}, 1/3, 1)$ we have
$$\dim\,{\HH}^2(\L) = \left\{
\begin{array}{ll}
1 & \mbox{ if $m\geq 3$ and $\kar K \neq 2$},\\
3 & \mbox{ if $m\geq 3$ and $\kar K = 2$,}\\
2 & \mbox{ if $m=2$ and $\kar K \neq 2$},\\
4 & \mbox{ if $m=2$ and $\kar K = 2$.}\\
\end{array}\right.
$$
\end{thm}

\begin{proof}
We consider first the case $m \geq 3$. Keeping the notation of
\ref{bit6} and Proposition \ref{D_{3m},1}, the set $f^3$ may be
chosen to consist of the following elements:
$$\{f^3_1, f^3_t, f^3_{m-1}, f^3_m\} \mbox{ with } t \in \{2, \ldots, m -
2\} \mbox{  where }$$
$$\begin{array}{l c l c l}
f^3_1 & = & f^2_1 \beta \alpha_1 \alpha_2 & = & \beta f^2_1 \alpha_1
\alpha_2 + \beta \alpha_1 \cdots \alpha_{m-1} f^2_2 \alpha_2 - \alpha_1
f^2_{3,2} \in e_1 K{\mathcal Q} e_3,\\

f^3_t & = & f^2_{3,t} \alpha_{t+1} & = & \alpha_t f^2_{3,t+1} \in e_t
K{\mathcal Q} e_{t+1} \mbox{ for }t \in \{2, \ldots, m-2\}\\

f^3_{m-1} & = & f^2_{3,m-1} \alpha_m & = & \alpha_{m-1} f^2_2 \alpha_2
\cdots \alpha_m \beta + \alpha_{m-1} \alpha_m f^2_1 \beta -  \alpha_{m-1}
\alpha_m \beta f^2_1 \in e_{m-1} K{\mathcal Q} e_1,\\

f^3_m & = & f^2_2 \alpha_2 \cdots \alpha_m \beta \alpha_1 & = &  - \alpha_m
f^2_1 \beta \alpha_1 + \alpha_m \beta f^2_1 \alpha_1 + \alpha_m \beta
\alpha_1 \cdots \alpha_{m-1} f^2_2 \in e_m K{\mathcal Q} e_2.\\
\end{array}$$

To find ${\Im}\,d_2$, let $f \in {\Hom}(Q^1, \L)$ and so
$$f(e_1 \otimes_{\beta} e_1) = c_1 e_1 + c_2 \beta + c_3 \beta^2 + c_4
\beta^3,$$
$$f(e_{1} \otimes_{\alpha_{1}} e_{2}) = d_{1} \alpha_{1} + k_1\beta\alpha_1,$$
$$f(e_{l} \otimes_{\alpha_{l}} e_{l +1}) = d_{l} \alpha_{l}, \mbox{ for } l
\in \{2, \ldots, m-1\},$$
$$f(e_m \otimes_{\alpha_m} e_1) = d_m \alpha_m + k_m\alpha_m\beta,$$
where $c_1, c_2, c_3, c_4, d_{l}, k_1, k_m \in K \mbox{ for } l \in
\{1, \ldots, m\}.$

It is straightforward to show that $fA_2$ is given by
$$fA_2(e_1 \otimes_{f^2_1} e_1) = 2c_1 \beta - (d_1 + d_2 + \ldots + d_m - 2
c_2) \beta^2 + (2c_3-k_1-k_m)\beta^3,$$
$$fA_2(e_m \otimes_{f^2_2} e_2) = (k_1+k_m)\alpha_m\beta\alpha_1,$$
$$f(e_j \otimes_{f^2_{3,j}} e_{j+1}) = 0, \mbox{ for
all } j \in \{2, \ldots, m-1\}.$$
So
$$\dim\,{\Im}\,d_2 = \left\{
\begin{array}{ll}
4 & \mbox{ if } \kar K \neq 2,\\
2 & \mbox{ if } \kar K = 2.\\
\end{array}\right.
$$

Now let $h \in {\Ker}\,d_3$, so  $h \in {\Hom}(Q^2, \L)$ and $d_3h = 0$.
Then $h: Q^2 \rightarrow \L$ is given by
$$h(e_1 \otimes_{f_1^2} e_1) = c_1 e_1 + c_2 \beta + c_3 \beta^2 + c_4
\beta^3,$$
$$h(e_m \otimes_{f_2^2} e_2) = c_5 \alpha_m \beta \alpha_1 \mbox{ and }$$
$$h(e_j \otimes_{f^2_{3,j}} e_{j+1}) = d_j \alpha_j, \mbox{ for } j \in \{2,
\ldots, m-1\},$$
for some $c_1, \ldots, c_5, d_j \in K \mbox{ where } j = 2, \ldots, m-1.$

By considering $hA_3(e_1 \otimes_{f^3_1} e_3)$ we see that $d_2 = 0.$

Then, for $t \in \{2, \ldots, m - 2\},$ we have $hA_3(e_t \otimes_{f^3_t}
e_{t+2}) = (d_t - d_{t+1}) \alpha_t \alpha_{t+1}.$ Then $d_t - d_{t+1} = 0$
and so $d_t = d_{t+1}$ for $t = 2, \ldots, m-2.$ Hence $d_2 = d_3 = \ldots =
d_{m-2} = d_{m-1}$. We already have $d_2 = 0$ so $d_j = 0$ for $j = 2,
\ldots, m-1.$

Moreover, $hA_3(e_m \otimes_{f^3_m} e_{2}) = 0$ so this gives us no
information. Thus, it may be verified that $h \in {\Ker}\,d_3$ is
given by
$$h(e_1 \otimes_{f_1^2} e_1) = c_1 e_1 + c_2 \beta + c_3 \beta^2 + c_4
\beta^3,$$
$$h(e_m \otimes_{f_2^2} e_2) = c_5 \alpha_m \beta \alpha_1 \mbox{ and }$$
$$h(e_j \otimes_{f^2_{3,j}} e_{j+1}) = 0, \mbox{ for } j \in \{2, \ldots,
m-1\}$$
for some $c_1, \ldots, c_5 \in K$ and so $\dim\,{\Ker}\,d_3 = 5.$

Therefore,
$$\dim\,{\HH}^2(\L) = \left\{
\begin{array}{ll}
5 - 4 = 1 & \mbox{ if } \kar K \neq 2,\\
5 - 2 = 3 & \mbox{ if } \kar K = 2.\\
\end{array}\right.
$$

For $m=2$, we again have that
$$\dim\,{\Im}\,d_2 = \left\{
\begin{array}{ll}
4 & \mbox{ if } \kar K \neq 2,\\
2 & \mbox{ if } \kar K = 2.\\
\end{array}\right.
$$
However, in this case we have that $\dim\,\Ker\, d_3 = 6$. Hence,
for $m=2$, we have
$$\dim\,{\HH}^2(\L) = \left\{
\begin{array}{ll}
2 & \mbox{ if } \kar K \neq 2,\\
4 & \mbox{ if } \kar K = 2.\\
\end{array}\right.
$$
This completes the proof.
\end{proof}

\vspace{1cm}

\begin{bit}\label{basisD_{3m}}
{\it A basis for ${\HH}^2(\L)$ for the standard algebra $\L =
\L(D_{3m}, 1/3,1)$ for $m \geq 3$.}
\end{bit}

\noindent {\bf Suppose $\kar K \neq 2.$}

From Theorem \ref{sumch8} we know that $\dim\,{\HH}^2(\L) = 1$ in
this case.  Let $h$ be the map given by
$$\begin{array}{rcl}
e_1 \otimes_{f_1^2} e_1 & \mapsto & e_1,\\
\mbox{else } & \mapsto & 0.\\
\end{array}$$
Then $\{h + {\Im}\,d_2\}$ is a basis of ${\HH}^2(\L)$ when $\kar K
\neq 2.$

\vspace{.5cm}

\noindent {\bf Suppose $\kar K = 2.$}

Here $\dim\,{\HH}^2(\L) = 3$ from Theorem \ref{sumch8}. We start by
defining non-zero maps $h_1, h_2, h_3$ in ${\Ker}\,d_3$.

Let $h_1$ be the map given by
$$\begin{array}{rcl}
e_1 \otimes_{f_1^2} e_1 & \mapsto & e_1,\\
\mbox{else } & \mapsto & 0,\\
\end{array}$$
$h_2$ be given by
$$\begin{array}{rcl}
e_1 \otimes_{f_1^2} e_1 & \mapsto & \beta,\\
\mbox{else } & \mapsto & 0,\\
\end{array}$$
and $h_3$ be given by
$$\begin{array}{rcl}
e_1 \otimes_{f_1^2} e_1 & \mapsto & \beta^3,\\
\mbox{else } & \mapsto & 0.\\
\end{array}$$

It can be shown that these maps are not in ${\Im}\,d_2$ since $\kar
K = 2$. Now we will show that $\{h_1 + {\Im}\,d_2, h_2 + {\Im}\,d_2,
h_3 + {\Im}\,d_2\}$ is a linearly independent set in ${\Ker}\,d_3 /
{\Im}\,d_2 = {\HH}^2(\L).$

Suppose $a(h_1 + {\Im}\,d_2) + b(h_2 + {\Im}\,d_2) + c(h_3 +
{\Im}\,d_2) = 0 + {\Im}\,d_2$ for some $a, b, c \in K$. So $a h_1 +
b h_2 + c h_3 \in {\Im}\,d_2$. Hence $a h_1 + b h_2 + c h_3 = fA_2$
for some $f \in {\Hom}(Q^1, \L)$.

Then $(a h_1 + b h_2 + c h_3)(e_1 \otimes_{f_1^2} e_1) = fA_2(e_1
\otimes_{f_1^2} e_1).$ So $a e_1 + b \beta + c \beta^3 = d\beta^2 -
k\beta^3$ for some $d, k \in K$. Since $\{e_1, \beta, \beta^2,
\beta^3 \}$ is linearly independent in $\L$, we have $a = b = 0$ and
$c = k.$ But $0 = (a h_1 + b h_2 + c h_3)(e_m \otimes_{f_2^2} e_2) =
fA_2(e_m \otimes_{f_2^2} e_2) = k\alpha_m\beta\alpha_1.$ So $k=0$
and thus $c=0$. Hence $\{h_1 + {\Im}\,d_2, h_2 + {\Im}\,d_2, h_3 +
{\Im}\,d_2\}$ is linearly independent in ${\HH}^2(\L)$ and forms a
basis of ${\HH}^2(\L)$ when $\kar K = 2.$

\vspace{.5cm}

\begin{bit}\label{basisD_{3m}m=2}
{\it A basis for ${\HH}^2(\L)$ for the standard algebra $\L =
\L(D_{3m}, 1/3,1)$ for $m = 2$.}
\end{bit}

Note first that $f^2_1 = \beta^2-\alpha_1\alpha_2$ and $f^2_2 =
\alpha_2\alpha_1$.

\noindent {\bf Suppose $\kar K \neq 2.$}

From Theorem \ref{sumch8} we know that $\dim\,{\HH}^2(\L) = 2$ in
this case.  Let $h_1$ be the map given by
$$\begin{array}{rcl}
e_1 \otimes_{f_1^2} e_1 & \mapsto & e_1,\\
\mbox{else } & \mapsto & 0,\\
\end{array}$$
and $h_2$ be given by
$$\begin{array}{rcl}
e_2 \otimes_{f_2^2} e_2 & \mapsto & e_2,\\
\mbox{else } & \mapsto & 0.\\
\end{array}$$
A similar argument to that above shows that $\{h_1 + {\Im}\,d_2, h_2
+ {\Im}\,d_2\}$ is a basis of ${\HH}^2(\L)$ when $\kar K \neq 2.$

\vspace{.5cm}

\noindent {\bf Suppose $\kar K = 2.$}

Here $\dim\,{\HH}^2(\L) = 4$ from Theorem \ref{sumch8}. Let $h_1$ be
the map given by
$$\begin{array}{rcl}
e_1 \otimes_{f_1^2} e_1 & \mapsto & e_1,\\
\mbox{else } & \mapsto & 0,\\
\end{array}$$
$h_2$ be given by
$$\begin{array}{rcl}
e_1 \otimes_{f_1^2} e_1 & \mapsto & \beta,\\
\mbox{else } & \mapsto & 0,\\
\end{array}$$
$h_3$ be given by
$$\begin{array}{rcl}
e_1 \otimes_{f_1^2} e_1 & \mapsto & \beta^3,\\
\mbox{else } & \mapsto & 0,\\
\end{array}$$
and $h_4$ be given by
$$\begin{array}{rcl}
e_2 \otimes_{f_2^2} e_2 & \mapsto & e_2,\\
\mbox{else } & \mapsto & 0.\\
\end{array}$$
Again, a similar argument shows that $\{h_1 + {\Im}\,d_2, h_2 +
{\Im}\,d_2, h_3 + {\Im}\,d_2, h_4 + {\Im}\,d_2\}$ is linearly
independent in ${\HH}^2(\L)$ and forms a basis of ${\HH}^2(\L)$ when
$\kar K = 2.$

\vspace{1cm}

\begin{thm}\label{sumch9}
For $\L = \L(E_n, 1, 1)$ with $n = 6, 7, 8$, we have ${\HH}^2(\L) = 0.$
\end{thm}

\begin{proof}
For $\L = \L(E_n, 1, 1)$ we have the quiver ${\mathcal Q}(E_n, 1)$ which is
described:
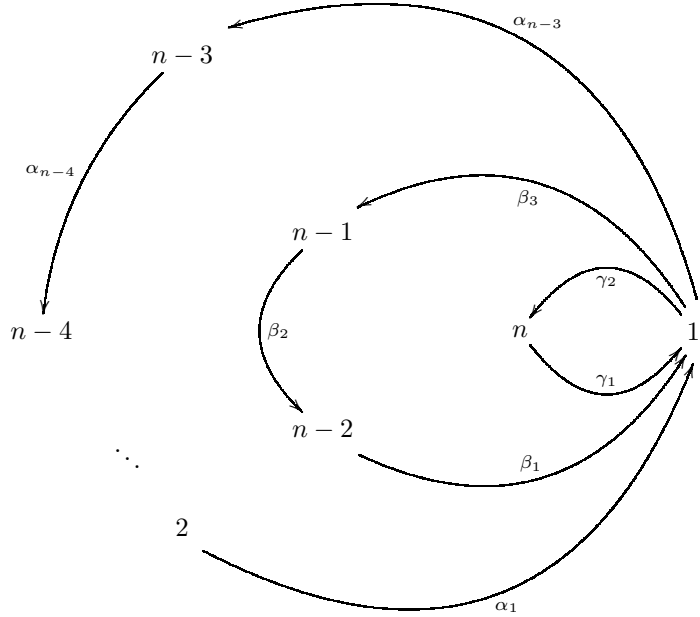
\vspace{1cm}

$$\xymatrix{
& & & & & & & \\
& {n-3}\ar@/_1pc/[dddl]_{\alpha_{n-4}} & & & & & & \\
& & & & & & & \\
&  & {n-1}\ar@/_2pc/[dd]^{\beta_2} & & & & & \\
{n-4} & & & & n \ar@/_2pc/[rr]^{\gamma_1} & &
1\ar@/_2pc/[ll]^{\gamma_2}\ar@<-.7ex>@/_3pc/[ullll]^{\beta_3}\ar@<-1.5ex>@/_5pc/[uuulllll]_{\alpha_{n-3}}
& \\
\ar@{}[dr]^\ddots & & {n-2}\ar@<-.7ex>@/_3pc/[urrrr]^{\beta_1} & & & & &  \\
& 2 \ar@<-1.5ex>@/_5pc/[uurrrrr]_{\alpha_1} & & & & & & \\
}$$

\vspace{1cm}

The set $f^3$ may be chosen to consist of the following elements:
$$\{f^3_{1,1}, f^3_{1,2}, f^3_{1,3}, f^3_{1,4}, f^3_{1,5}, f^3_2, f^3_{n-1},
f^3_{n-2}, f^3_n, f^3_3, f^3_m\} \mbox{ where }$$
$$\begin{array}{l c l c l}
f^3_{1,1} & = & f^2_{1,2} \alpha_{n-3} \alpha_{n-4} & = & \beta_3 \beta_2
f^2_{2,3} \alpha_{n-4} - \alpha_{n-3} f^2_{3,n-4} \in e_1 K{\mathcal Q}
e_{n-4},\\
f^3_{1,2} & = & f^2_{1,1} \alpha_{n-3} & = & \beta_3 \beta_2 f^2_{2,3} -
\gamma_2 f^2_{2,5} \in e_1 K{\mathcal Q} e_{n-3},\\
f^3_{1,3} & = & f^2_{1,1} \beta_3\beta_2 & = & \beta_3 f^2_4 - \gamma_2
f^2_{2,6} \beta_2 \in e_1 K{\mathcal Q} e_{n-2},\\
f^3_{1,4} & = & f^2_{1,1} \beta_3 - f^2_{1,2} \beta_3 & = &
\alpha_{n-3}\alpha_{n-4} \cdots \alpha_2 f^2_{2,1} - \gamma_2 f^2_{2,6} \in
e_1 K{\mathcal Q} e_{n-1},\\
f^3_{1,5} & = & f^2_{1,2} \gamma_2 & = & \beta_3\beta_2 f^2_{2,4} -
\alpha_{n-3} \alpha_{n-4} \cdots \alpha_2 f^2_{2,2} \in e_1 K{\mathcal Q}
e_n,\\
f^3_2 & = & f^2_{2,1} \beta_2\beta_1 - f^2_{2,2} \gamma_1 & = & \alpha_1
f^2_{1,1} \in e_2 K{\mathcal Q} e_1,\\
f^3_{n-1} & = & f^2_4 \beta_1 & = & \beta_2\beta_1 f^2_{1,1} +  \beta_2
f^2_{2,4}\gamma_1 \in e_{n-1} K{\mathcal Q} e_1,\\
f^3_{n-2} & = & f^2_{2,3} \alpha_{n-4} \cdots \alpha_2 \alpha_1 - f^2_{2,4}
\gamma_1 & = & \beta_1 f^2_{1,1} - \beta_1 f^2_{1,2} \in e_{n-2} K{\mathcal
Q} e_1,\\
f^3_n & = & f^2_{2,6} \beta_2 \beta_1 - f^2_{2,5} \alpha_{n-4} \cdots
\alpha_2\alpha_1 & = & \gamma_1 f^2_{1,2} \in e_n K{\mathcal Q} e_1,\\
f^3_3 & = & f^2_{3,2} \alpha_1 & = & \alpha_2 f^2_{2,1} \beta_2 \beta_1 -
\alpha_2 \alpha_1 f^2_{1,2} \in e_3 K{\mathcal Q} e_1,\\
f^3_m & = & f^2_{3,m-1} \alpha_{m-2} & = & \alpha_{m-1} f^2_{3,m-2}  \in e_m
K{\mathcal Q} e_{m-2},\\
&&&&\mbox{ for } m = 4, \ldots n-3.
\end{array}$$
Then it is easy to check by direct calculations that $\dim\,{\Ker}\,d_3 =
\dim\,{\Im}\,d_2 = 2$ and so ${\HH}^2(\L) = 0$.
\end{proof}

\begin{thm}\label{sumch10}
For $\L = \L(E_6, 1, 2)$ we have ${\HH}^2(\L) = 0.$
\end{thm}

\begin{proof}
With the notation for ${\mathcal Q}(E_6, 1)$ as in Theorem \ref{sumch9} and
with $n = 6$, the set $f^3$ may be chosen to consist of the following
elements:
$$\{f^3_{1,1}, f^3_{1,2}, f^3_{1,3}, f^3_{1,4}, f^3_{1,5}, f^3_2, f^3_3,
f^3_4, f^3_5, f^3_6\} \mbox{ where }$$
$$\begin{array}{l c l c l}
f^3_{1,1} & = & f^2_{1,2} \gamma_2 & = & \beta_3 \beta_2 f^2_{2,4} -
\alpha_3 \alpha_2 f^2_{2,3} \in e_1 K{\mathcal Q} e_6,\\
f^3_{1,2} & = & f^2_{1,1} \beta_3 & = &  \beta_3\beta_2 f^2_{2,6} - \gamma_2
f^2_{2,2} \in e_1 K{\mathcal Q} e_5,\\
f^3_{1,3} & = & f^2_{1,2} \beta_3 \beta_2 - \beta_3 \beta_2 f^2_{2,6}
\beta_2 & = & - \alpha_3 f^2_{3,1} \in e_1 K{\mathcal Q} e_4,\\
f^3_{1,4} & = & f^2_{1,1} \alpha_3 - f^2_{1,2} \alpha_3 & = & \alpha_3
\alpha_2 f^2_{2,5} - \gamma_2 f^2_{2,1} \in e_1 K{\mathcal Q} e_3,\\
f^3_{1,5} & = & f^2_{1,1} \alpha_3 \alpha_2 + \gamma_2 f^2_{2,1} \alpha_2 &
= & \beta_3 f^2_{3,2} \in e_1 K{\mathcal Q} e_2,\\
f^3_2 & = & f^2_{2,5} \alpha_2 \alpha_1 - f^2_{2,3} \gamma_1 & = & \alpha_1
f^2_{1,1} - \alpha_1 f^2_{1,2} \in e_2 K{\mathcal Q} e_1,\\
f^3_3 & = & f^2_{3,1} \beta_1 - \alpha_2 f^2_{2,3} \gamma_1 & = & \alpha_2
\alpha_1 f^2_{1,1} \in e_3 K{\mathcal Q} e_1,\\
f^3_4 & = & f^2_{2,6} \beta_2 \beta_1 - f^2_{2,4} \gamma_1 & = & \beta_1
f^2_{1,1} \in e_4 K{\mathcal Q} e_1,\\
f^3_5 & = & f^2_{3,2} \alpha_1 - \beta_2 f^2_{2,6} \beta_2 \beta_1 & = & -
\beta_2 \beta_1 f^2_{1,2} \in e_5 K{\mathcal Q} e_1,\\
f^3_6 & = & f^2_{2,2} \beta_2\beta_1 - f^2_{2,1} \alpha_2 \alpha_1 & = &
\gamma_1 f^2_{1,2} \in e_6 K{\mathcal Q} e_1.
\end{array}$$
Again by direct calculations we can show that $\dim\,{\Ker}\,d_3 =
\dim\,{\Im}\,d_2 = 2$ and so ${\HH}^2(\L) = 0$.
\end{proof}

\vspace{1cm}

To summarise the results of Sections \ref{sec4} and \ref{sec5} we have the
following theorem.

\begin{thm}\label{sumthm}
Let $\L$ be a standard self-injective algebra of finite representation type
of type $\L(D_n, s, 1), \L(D_4, s, 3)$ with $n \geq 4, s \geq 1$, $\L(D_n,
s, 2), \L(D_{3m}, s/3, 1)$, $3 \nmid s$ with $n \geq 4, m \geq 2, s \geq 2$
or $\L(E_n, s, 1)$, $\L(E_6, s, 2)$ with $n \in \{6, 7, 8\}, s \geq 1$. Then
${\HH}^2(\L) = 0.$

Let $\L$ be $\L(D_n, 1, 2)$; then $\dim\,{\HH}^2(\L) = 1$ and a basis for
${\HH}^2(\L)$ is given in \ref{basisD_n}.

Let $\L$ be $\L(D_{3m}, 1/3, 1)$; then
$$\dim\,{\HH}^2(\L) = \left\{
\begin{array}{ll}
1 & \mbox{ if $m \geq 3$ and $\kar K \neq 2$,}\\
3 & \mbox{ if $m \geq 3$ and $\kar K = 2,$}\\
2 & \mbox{ if $m = 2$ and $\kar K \neq 2,$}\\
4 & \mbox{ if $m = 2$ and $\kar K = 2.$}\\
\end{array}\right.
$$
and a basis for ${\HH}^2(\L)$ is given in \ref{basisD_{3m}} and
\ref{basisD_{3m}m=2}.
\end{thm}

\vspace{1.5cm}

Thus with the information taken from \cite{EH} and \cite{GS} for the
algebras of type $A_n$, we now know the second Hochschild cohomology group
for all standard finite dimensional self-injective algebras of finite
representation type over an algebraically closed field $K$.

\section{${\HH}^2(\L)$ for the non-standard self-injective algebras of
finite representation type}\label{sec6}

Let $\L = \L(m)$, $m \geq 2$, be the non-standard algebra of
\ref{bit9} so we assume now that the characteristic of $K$ is 2. We
may choose a minimal generating set $f^2$ with elements as follows:
$$f_1^2 = \beta^2 - \alpha_1 \cdots \alpha_m, \hspace{1cm} f_2^2 = \alpha_m
\alpha_1 - \alpha_m \beta \alpha_1,$$
$$f_{3,j}^2 = \alpha_j \alpha_{j+1} \cdots \alpha_j
\mbox { for }  \left\{
\begin{array}{ll}
j = 2, \ldots, m-1  & \mbox{ if } m \geq 3,\\
j = 2  & \mbox{ if } m = 2.\\
\end{array}\right.
$$

We know that ${\HH}^2(\L) = {\Ker}\,d_3 / {\Im}\,d_2$. First we will find
${\Im}\,d_2$. Let $f \in {\Hom}(Q^1, \L)$ and so
$$f(e_1 \otimes_{\beta} e_1) = c_1 e_1 + c_2 \beta + c_3 \beta^2 + c_4
\beta^3,$$
$$f(e_1 \otimes_{\alpha_1} e_{2}) = d_1 \alpha_1 + k_1\beta\alpha_1, $$
$$f(e_l \otimes_{\alpha_l} e_{l +1}) = d_l \alpha_l, \mbox{ for } l \in \{2,
\ldots, m-1\},$$
$$f(e_m \otimes_{\alpha_m} e_1) = d_m \alpha_m + k_m\alpha_m\beta,$$
where $c_1, c_2, c_3, c_4, d_l, k_1, k_m \in K$ for $l \in \{1,
\ldots, m\}.$

We have
$Q^2 = (\L e_1 \otimes_{f^2_1} e_1\L) \oplus  (\L e_m \otimes_{f^2_2} e_2
\L) \oplus \bigoplus_{j = 2}^{m - 1} (\L e_j \otimes_{f^2_{3,j}} e_{j+1}
\L)$ if $m \geq 3$ and $Q^2 = (\L e_1 \otimes_{f^2_1} e_1\L) \oplus  (\L e_2
\otimes_{f^2_2} e_2 \L) \oplus (\L e_2 \otimes_{f^2_{3,2}} e_3 \L)$ if $m =
2.$

Now we find $fA_2$. We have $fA_2(e_1 \otimes_{f^2_1} e_1) = f(e_1
\otimes_{\beta} e_1) \beta + \beta f(e_1 \otimes_{\beta} e_1) -
f(e_1 \otimes_{\alpha_1} e_2) \alpha_2 \cdots \alpha_m - \alpha_1
f(e_2 \otimes_{\alpha_2} e_3) \alpha_3 \cdots \alpha_m - \cdots -
\alpha_1 \alpha_2 \cdots \alpha_{m-1} f(e_m \otimes_{\alpha_m} e_1)
= (c_1 e_1 + c_2 \beta + c_3 \beta^2 + c_4 \beta^3) \beta +
\beta(c_1 e_1 + c_2 \beta + c_3 \beta^2 + c_4 \beta^3) - d_1
\alpha_1 \cdots \alpha_m - d_2 \alpha_1 \cdots \alpha_m - \ldots -
d_m \alpha_1 \cdots \alpha_m - k_1\beta\alpha_1 \cdots \alpha_m
-k_m\alpha_1 \cdots \alpha_m\beta = 2 c_1 \beta - (d_1 + d_2 +
\cdots + d_m - 2c_2) \beta^2 + (2c_3 - k_1 - k_m)\beta^3.$

Also $fA_2(e_m \otimes_{f^2_2} e_2) = f(e_m \otimes_{\alpha_m} e_1)
\alpha_1 + \alpha_m f(e_1 \otimes_{\alpha_1} e_2) - f(e_m
\otimes_{\alpha_m} e_1) \beta \alpha_1 - \alpha_m f(e_1
\otimes_{\beta} e_1) \alpha_1 - \alpha_m \beta f(e_1
\otimes_{\alpha_1} e_2) = (d_m \alpha_m + k_m\alpha_m\beta)\alpha_1
+ \alpha_m(d_1 \alpha_1 + k_1\beta\alpha_1) - (d_m \alpha_m +
k_m\alpha_m\beta)\beta \alpha_1 - \alpha_m(c_1 e_1 + c_2 \beta + c_3
\beta^2 + c_4 \beta^3) \alpha_1 - \alpha_m \beta (d_1\alpha_1 +
k_1\beta\alpha_1) = (k_1 + k_m - c_1 - c_2) \alpha_m \alpha_1.$

Finally, for $m \geq 3$ and $j = 2, \ldots, m-1$ or for $m = 2$ and
$j = 2$, we have $fA_2(e_j \otimes_{f^2_{3,j}} e_{j+1}) = 0.$

Thus $fA_2$ is given by
$$fA_2(e_1 \otimes_{f^2_1} e_1) = 2c_1 \beta - (d_1 + d_2 + \cdots + d_m - 2
c_2) \beta^2 + (2c_3 - k_1 - k_m) \beta^3,$$
$$fA_2(e_m \otimes_{f^2_2} e_2) = (k_1 + k_m - c_1 - c_2) \alpha_m \alpha_1,$$
$$f(e_j \otimes_{f^2_{3,j}} e_{j+1}) = 0, \mbox{ for all } j \in \{2,
\ldots, m-1\} \mbox{ if } m \geq 3 \mbox{ or } j = 2 \mbox{ if } m =
2.$$ So, since $\kar K = 2$, we have $\dim\,{\Im}\,d_2 = 3$.

\vspace{1cm}

Next we determine ${\Ker}\,d_3$. We need to consider separately the
cases $m \geq 3$ and $m = 2$. Suppose first that $m \geq 3$.

\vspace{.5cm}

For $m \geq 3$, we choose the set $f^3$ to consist of the following
elements:
$$\{f^3_1, f^3_t, f^3_{m-1}, f^3_m\} \mbox{ with } t \in \{2, \ldots, m -
2\} \mbox{  where }$$
$$\begin{array}{l c l c l}
f^3_1 & = & f^2_1 \beta \alpha_1 \alpha_2 & = & \beta f^2_1 \alpha_1
\alpha_2 + \alpha_1 \cdots \alpha_{m-1} f^2_2 \alpha_2 + (\beta \alpha_1 -
\alpha_1) f^2_{3,2} \in e_1 K{\mathcal Q} e_3,\\

f^3_t & = & f^2_{3,t} \alpha_{t+1} & = & \alpha_t f^2_{3,t+1} \in e_t
K{\mathcal Q} e_{t+2} \mbox { for } t \in \{2, \ldots, m-2\},\\

f^3_{m-1} & = & f^2_{3,m-1} (\alpha_m - \alpha_m \beta) & = & \alpha_{m-1}
f^2_2 \alpha_2 \cdots \alpha_m  + \alpha_{m-1} \alpha_m f^2_1 \beta -
\alpha_{m-1} \alpha_m \beta f^2_1 \in e_{m-1} K{\mathcal Q} e_1,\\

f^3_m & = & f^2_2 \alpha_2 \cdots \alpha_m \alpha_1 & = &  - \alpha_m f^2_1
\beta \alpha_1 + \alpha_m \beta f^2_1 \alpha_1 + \alpha_m \alpha_1 \cdots
\alpha_{m-1} f^2_2 \in e_m K{\mathcal Q} e_2.\\
\end{array}$$

\vspace{.5cm}

Let $h \in {\Ker}\,d_3$. Then $h: Q^2 \rightarrow \L$ is given by
$$h(e_1 \otimes_{f_1^2} e_1) = c_1 e_1 + c_2 \beta + c_3 \beta^2 + c_4
\beta^3,$$
$$h(e_m \otimes_{f_2^2} e_2) = c_5 \alpha_m \alpha_1 \mbox{ and }$$
$$h(e_j \otimes_{f^2_{3,j}} e_{j+1}) = d_j \alpha_j, \mbox{ for } j \in \{2,
\ldots, m-1\},$$
for some $c_1, \ldots, c_5, d_j \in K \mbox{ where } j = 2, \ldots, m-1.$

Then $hA_3(e_1 \otimes_{f^3_1} e_3) = h(e_1 \otimes_{f^2_1} e_1) \beta
\alpha_1 \alpha_2 - \beta h(e_1 \otimes_{f^2_1} e_1) \alpha_1 \alpha_2 -
\alpha_1 \cdots \alpha_{m-1} h(e_m \otimes_{f^2_2} e_2)\alpha_2 - (\beta
\alpha_1 - \alpha_1) h(e_2 \otimes_{f^2_{3,2}} e_3) = (c_1 e_1 + c_2 \beta +
c_3 \beta^2 + c_4 \beta^3) \beta \alpha_1 \alpha_2 - \beta (c_1 e_1 + c_2
\beta + c_3 \beta^2 + c_4 \beta^3) \alpha_1 \alpha_2 - c_5 \alpha_1 \cdots
\alpha_{m-1} \alpha_m \alpha_1 \alpha_2 - d_2 \beta \alpha_1 \alpha_2 + d_2
\alpha_1 \alpha_2 = d_2 (\alpha_1 \alpha_2 - \beta \alpha_1 \alpha_2).$ As
$h \in {\Ker}\,d_3$ we have $d_2 = 0.$

For $t \in \{2, \ldots, m - 2\},$ we have $hA_3(e_t \otimes_{f^3_t} e_{t+2})
= h(e_t \otimes_{f^2_{3,t}} e_{t+1}) \alpha_{t+1} - \alpha_t h(e_{t+1}
\otimes_{f^2_{3, t+1}} e_{t+2}) = d_t \alpha_t \alpha_{t+1} - d_{t+1}
\alpha_t \alpha_{t+1} = (d_t - d_{t+1}) \alpha_t \alpha_{t+1}.$ Then $d_t -
d_{t+1} = 0$ and so $d_t = d_{t+1}$ for $t = 2, \ldots, m-2.$ Hence $d_2 =
d_3 = \ldots = d_{m-2} = d_{m-1}$. We already have $d_2 = 0$ so $d_j = 0$
for $j = 2, \ldots, m-1.$

Now $hA_3(e_{m-1} \otimes_{f^3_{m-1}} e_1) = h(e_{m-1} \otimes_{f^2_{3,m-1}}
e_m) (\alpha_m - \alpha_m \beta) - \alpha_{m-1} h(e_m \otimes_{f^2_2}
e_2)\alpha_2 \cdots \alpha_m - \alpha_{m-1} \alpha_m h(e_1 \otimes_{f^2_1}
e_1)\beta + \alpha_{m-1} \alpha_m \beta h(e_1 \otimes_{f^2_1} e_1) = d_{m-1}
\alpha_{m-1} \alpha_m - d_{m-1} \alpha_{m-1} \alpha_m \beta - c_5
\alpha_{m-1} \alpha_m \alpha_1 \alpha_2 \cdots \alpha_m - \alpha_{m-1}
\alpha_m (c_1 e_1 + c_2 \beta + c_3 \beta^2 + c_4 \beta^3) \beta +
\alpha_{m-1} \alpha_m \beta (c_1 e_1 + c_2 \beta + c_3 \beta^2 + c_4
\beta^3) = d_{m-1}(\alpha_{m-1} \alpha_m - \alpha_{m-1} \alpha_m \beta) =
0,$ as $d_{m-1} = 0$ from above.

Finally, $hA_3(e_m \otimes_{f^3_m} e_2) = h(e_m \otimes_{f^2_2} e_2)
\alpha_2 \cdots \alpha_m \alpha_1 + \alpha_m h(e_1 \otimes_{f^2_1}
e_1) \beta \alpha_1 - \alpha_m \beta h(e_1 \otimes_{f^2_1}
e_1)\alpha_1 - \alpha_m \alpha_1 \cdots \alpha_{m-1} h(e_m
\otimes_{f^2_2} e_2) = c_5 \alpha_m \alpha_1 \alpha_2 \cdots
\alpha_m \alpha_1 + \alpha_m (c_1 e_1 + c_2 \beta + c_3 \beta^2 +
c_4 \beta^3) \beta \alpha_1 - \alpha_m \beta (c_1 e_1 + c_2 \beta +
c_3 \beta^2 + c_4 \beta^3)\alpha_1 - c_5 \alpha_m \alpha_1 \cdots
\alpha_{m-1} \alpha_m \alpha_1 = -c_1 \alpha_m \beta \alpha_1 + c_1
\alpha_m \beta \alpha_1 = 0,$ and so this gives no information on
the constants occurring in $h$.

Thus $h$ is given by
$$h(e_1 \otimes_{f_1^2} e_1) = c_1 e_1 + c_2 \beta + c_3 \beta^2 + c_4
\beta^3,$$
$$h(e_m \otimes_{f_2^2} e_2) = c_5 \alpha_m \alpha_1 \mbox{ and }$$
$$h(e_j \otimes_{f^2_{3,j}} e_{j+1}) = 0, \mbox{ for } j \in \{2, \ldots,
m-1\}$$
for some $c_1, \ldots, c_5 \in K$ and so $\dim\,{\Ker}\,d_3 = 5.$

Therefore, for $m \geq 3$ we have $\dim\,{\HH}^2(\L) = 5 - 3 = 2$.

\smallskip

This gives the following theorem.

\begin{thm}\label{sumch12}
For $\L = \L(m)$ and $m \geq 3$ we have $\dim\,{\HH}^2(\L) = 2.$
\end{thm}

\vspace{.5cm}

\begin{bit}\label{nonstandard basis}
{\it A basis for ${\HH}^2(\L)$ for $\L = \L(m)$ and $m \geq 3$.}
\end{bit}

\noindent We have $\kar K = 2$, $m \geq 3$, and $\dim\,{\HH}^2(\L) =
2$. We start by defining non-zero maps $h_1, h_2$ in ${\Ker}\,d_3$.

Let $h_1$ be the map given by
$$\begin{array}{rcl}
e_1 \otimes_{f_1^2} e_1 & \mapsto & e_1,\\
\mbox{else } & \mapsto & 0,\\
\end{array}$$
and $h_2$ be given by
$$\begin{array}{rcl}
e_1 \otimes_{f_1^2} e_1 & \mapsto & \beta,\\
\mbox{else } & \mapsto & 0.\\
\end{array}$$

It can be shown as before that these maps are not in ${\Im}\,d_2$.
Now we will show that $\{h_1 + {\Im}\,d_2, h_2 + {\Im}\,d_2\}$ is a
linearly independent set in ${\HH}^2(\L).$

Suppose $a(h_1 + {\Im}\,d_2) + b(h_2 + {\Im}\,d_2) = 0 + {\Im}\,d_2$
for some $a, b \in K$. So $a h_1 + b h_2 \in {\Im}\,d_2$. Hence $a
h_1 + b h_2 = fA_2$ for some $f \in {\Hom}(Q^1, \L)$. Then $(a h_1 +
b h_2)(e_1 \otimes_{f_1^2} e_1) = fA_2(e_1 \otimes_{f_1^2} e_1).$ So
$a e_1 + b \beta = d \beta^2 + k \beta^3$ for some $d, k \in K$.
Since $\{e_1, \beta, \beta^2, \beta^3 \}$ is linearly independent in
$\L$, we have $a = b = 0.$ Hence $\{h_1 + {\Im}\,d_2, h_2 +
{\Im}\,d_2\}$ is linearly independent in ${\HH}^2(\L)$ and forms a
basis of ${\HH}^2(\L)$.

\vspace{.5cm}

\begin{bit}\label{nonstandard basis m=2}
{\it ${\HH}^2(\L)$ in the case $\L = \L(m)$ and $m = 2$.}
\end{bit}

In the case $m = 2$ we showed above that $\dim\,{\Im}\,d_2 = 3$. But
now we have $\dim\,{\Ker}\,d_3 = 6$. Thus $\dim\,\HH^2(\L) = 3$. It
can be verified that $\{h_1 + {\Im}\,d_2, h_2 + {\Im}\,d_2, h_3 +
{\Im}\,d_2\}$ is a basis of ${\HH}^2(\L)$, where $h_1$ is the map
given by
$$\begin{array}{rcl}
e_1 \otimes_{f_1^2} e_1 & \mapsto & e_1,\\
\mbox{else } & \mapsto & 0,\\
\end{array}$$
$h_2$ is given by
$$\begin{array}{rcl}
e_1 \otimes_{f_1^2} e_1 & \mapsto & \beta,\\
\mbox{else } & \mapsto & 0.\\
\end{array}$$
and $h_3$ is given by
$$\begin{array}{rcl}
e_2 \otimes_{f_2^2} e_2 & \mapsto & e_2,\\
e_2 \otimes_{f_3^2} e_1 & \mapsto & \alpha_2 + \alpha_2\beta,\\
\mbox{else } & \mapsto & 0.\\
\end{array}$$

\vspace{1cm}

We summarise all these results in the following theorem.

\begin{thm}
For $\L = \L(m)$ where $\kar K = 2, m \geq 2$ we have
$$\dim\,{\HH}^2(\L) = \left\{
\begin{array}{ll}
2 & \mbox{ if $m \geq 3$,}\\
3 & \mbox{ if $m = 2$.}\\
\end{array}\right.$$

Moreover, if $m\geq 3$ then $\{h_1 + {\Im}\,d_2, h_2 + {\Im}\,d_2\}$
is a basis for ${\HH}^2(\L)$ where $h_1$ is the map given by
$$\begin{array}{rcl}
e_1 \otimes_{f_1^2} e_1 & \mapsto & e_1,\\
\mbox{else } & \mapsto & 0,\\
\end{array}$$
and $h_2$ is given by
$$\begin{array}{rcl}
e_1 \otimes_{f_1^2} e_1 & \mapsto & \beta,\\
\mbox{else } & \mapsto & 0.\\
\end{array}$$
If $m = 2$ then $\{h_1 + {\Im}\,d_2, h_2 + {\Im}\,d_2, h_3 +
{\Im}\,d_2\}$ is a basis for ${\HH}^2(\L)$ where $h_1$ is the map
given by
$$\begin{array}{rcl}
e_1 \otimes_{f_1^2} e_1 & \mapsto & e_1,\\
\mbox{else } & \mapsto & 0,\\
\end{array}$$
$h_2$ is given by
$$\begin{array}{rcl}
e_1 \otimes_{f_1^2} e_1 & \mapsto & \beta,\\
\mbox{else } & \mapsto & 0,\\
\end{array}$$
and $h_3$ is given by
$$\begin{array}{rcl}
e_2 \otimes_{f_2^2} e_2 & \mapsto & e_2,\\
e_2 \otimes_{f_3^2} e_1 & \mapsto & \alpha_2 + \alpha_2\beta,\\
\mbox{else } & \mapsto & 0.\\
\end{array}$$
\end{thm}

\vspace{1cm}

This completes the discussion of ${\HH}^2(\L)$ for the non-standard
self-injective algebras of finite representation type over an algebraically
closed field.

To conclude we now summarise ${\HH}^2(\L)$ for all finite dimensional
self-injective algebras of finite representation type over an algebraically
closed field.
\begin{thm}\label{finalthm}
Let $\L$ be a finite dimensional self-injective algebra of finite
representation type over an algebraically closed field $K$. If $\L$ is the
standard algebra of type $\L(A_{2p+1}, s, 2)$ with $s, p \geq 2$, $\L(D_n,
s, 1), \L(D_4, s, 3)$ with $n \geq 4, s \geq 1$, $\L(D_n, s, 2), \L(D_{3m},
s/3, 1)$ with $n \geq 4, m \geq 2, s \geq 2$ or $\L(E_n, s, 1)$, $\L(E_6, s,
2)$ with $n \in \{6, 7, 8\}, s \geq 1$; then ${\HH}^2(\L) = 0.$

If $\L$ is of type $\L(A_n, s/n, 1)$ then $\dim\,{\HH}^2(\L) = m$ where $n +
1 = ms + r$ and $0 \leq r < s$.

For $\L(A_3, 1, 2)$; then $\dim\,{\HH}^2(\L) = 1.$

Let $\L$ be $\L(D_n, 1, 2)$; then $\dim\,{\HH}^2(\L) = 1.$

Let $\L$ be the standard algebra $\L(D_{3m}, 1/3, 1)$; then
$$\dim\,{\HH}^2(\L) = \left \{ \begin{array}{ll}
1 & \mbox{ if $m \geq 3$ and $\kar K \neq 2$,}\\
3 & \mbox{ if $m \geq 3$ and $\kar K = 2,$}\\
2 & \mbox{ if $m = 2$ and $\kar K \neq 2,$}\\
4 & \mbox{ if $m = 2$ and $\kar K = 2.$}\\
\end{array} \right.$$

Let $\L$ be the non-standard algebra $\L(m)$ where $\kar K = 2, m
\geq 2$. Then $\dim\,{\HH}^2(\L) = 2$ if $m \geq 3$ and
$\dim\,{\HH}^2(\L) = 3$ if $m = 2$.
\end{thm}

\end{document}